\documentclass[a4paper,12pt,leqno,oneside]{amsart}

\usepackage{amssymb,hyperref}
\usepackage[latin1]{inputenc}
\swapnumbers
\theoremstyle{plain}
\newtheorem{maintheo}{Theorem}[section]
\newtheorem{btheo}[maintheo]{Theorem}
\newtheorem{ctheo}[maintheo]{Theorem}
\newtheorem{theo}{Theorem}[section]
\newtheorem{prop}[theo]{Proposition}
\newtheorem{cor}[theo]{Corollary}
\newtheorem{lem}[theo]{Lemma}

\newtheorem{sulem}[theo]{Sub-lemma}
\theoremstyle{definition}
\newtheorem{de}[theo]{Definition}
\newtheorem{notn}[theo]{Notation}
\newtheorem{claim}[theo]{Claim}
\theoremstyle{remark}
\newtheorem{rem}[theo]{Remark}

\numberwithin{equation}{section}

\newcommand{\thismonth}{\ifcase\month\or
  January\or February\or March\or April\or May\or June\or July\or
  August\or September\or October\or November\or December\fi
  \space\number\year}
\DeclareMathAlphabet{\mathrmsl}{OT1}{cmr}{m}{sl}

\newcommand{\oper}[3][n]{\newcommand{#2}{\mathop{\mathrm{#3}}%
\ifx n#1\nolimits\else\limits\fi} }
\newcommand{\rsoper}[3][n]{\newcommand{#2}{\mathop{\mathrmsl{#3}}%
\ifx n#1\nolimits\else\limits\fi} }
\newcommand{\ache}{asymptot­ically comp­lex hyper­bolic Eins­tein\ }
\newcommand{\chii}{{\mathbb C}\mathbf{H}^{2}}
\newcommand{\BM}{{\mathbb B}}
\newcommand{\Sym}{\operatorname{Sym}}
\newcommand{\NM}{{\mathbb N}}
\newcommand{\RM}{{\mathbb R}}
\newcommand{\CM}{{\mathbb C}}
\newcommand{\ra}{\rangle}
\newcommand{\la}{\langle}

\newcommand{\cJ}{\mathcal{J}}
\newcommand{\cQ}{\mathcal{Q}}
\newcommand{\cW}{\mathcal{W}}

\newcommand{\p}{\partial}
\newcommand{\pb}{\partial_{\infty}}
\newcommand{\vol}{\operatorname{vol}}
\newcommand{\Scal}{\operatorname{Scal}}
\newcommand{\Ric}{\operatorname{Ric}}
\newcommand{\tr}{\operatorname{tr}}
\newcommand{\Rsym}{\operatorname{\overset{\circ}{R}}}
\newcommand{\csum}{\mathfrak{S}}
\newcommand{\tv}{\mathcal{T}}

\newcommand{\ad}{\operatorname{ad}}
\newcommand{\II}{\mathbb{I}}
\renewcommand{\geq}{\geqslant}
\renewcommand{\leq}{\leqslant}
\newcommand{\cerc}{{\mathbb S}}
\newcommand{\eps}{\varepsilon}
\renewcommand{\exp}{\operatorname{e}}
\newcommand{\bg}{\bar g}
\newcommand{\dg}{\dot{g}}
\newcommand{\dR}{\dot{R}}
\newcommand{\cw}{C^{\infty}}  

\newcommand{\ti}{\theta^1}
\newcommand{\vti}{\vartheta^1}
\newcommand{\vt}{\vartheta^0}
\newcommand{\vtb}{\vartheta^{\bar{0}}}
\newcommand{\tib}{\theta^{\bar{1}}}
\newcommand{\vtib}{\vartheta^{\bar{1}}}
\newcommand{\oii}{\omega_1^1}
\newcommand{\taui}{\tau^1}
\newcommand{\tauib}{\tau^{\bar{1}}}
\newcommand{\dbar}{\overline{\partial}}
\newcommand{\wn}{\widetilde{\nabla}^W}
\newcommand{\proofof}[1]{\end{#1}\begin{proof}}

\newcounter{mnotecount}[section]
\renewcommand{\themnotecount}{\thesection.\arabic{mnotecount}}
\newcommand{\mnote}[1]
{\protect{\stepcounter{mnotecount}}$^{\mbox{\footnotesize  $
      \bullet$\themnotecount}}$ \marginpar{\raggedright\tiny\em
    $\!\!\!\!\!\!\,\bullet$\themnotecount: #1} }

\begin{document}

\title{A Burns-Epstein invariant for ACHE 4-manifolds}
\author{Olivier Biquard}
\address{IRMA, CNRS et Université Louis Pasteur\\
7 rue René Descartes\\ F-67084 Strasbourg Cedex\\France}
\author{Marc Herzlich}
\address{Département des Sciences Mathématiques\\
Géométrie--Topologie et Algèbre\\
UMR 5030 CNRS et Université Mont­pellier~II\\ 
Pl. E. Bataillon\\ F-34095 Montpellier Cedex 5\\France}

\begin{abstract}
We define a renormalized characteristic class for Einstein asymptotically
complex hyperbolic (\textsc{ache}) manifolds of dimension $4$: for any
such manifold, the polynomial in the curvature associated to the 
characteristic class $\chi -3\tau$ is shown to converge. This
extends a work of Burns and Epstein in the Kähler-Einstein case.

We also define a new global invariant for any compact 3-dimensional
pseudoconvex CR manifold, by a renormalization procedure of the $\eta$
invariant of a sequence of metrics which approximate the CR structure.

Finally, we get a formula relating the renormalized
characteristic class to the topological number $\chi -3\tau$ and the
invariant of the CR structure arising at infinity.
\end{abstract}


\thanks{Both authors are members of the {\sc edge} Research Training 
Network {\sc hprn-ct-2000-00101} of the European Union. The second author 
is also supported in part by an {\sc aci} program 
of the French Ministry of Research.}

\maketitle

\bigskip

\section{Introduction}

\bigskip

In \cite{BurEps90b}, Burns and Epstein showed
that, for complete Kähler-Einstein metrics on bounded domains in 
$\CM^m$ or in a complex manifold, the local integrands of some 
(precisely known) combinations of the Chern classes had convergent
integrals, thus providing interesting invariants of bounded domains.
In some cases, they were also able to compute {\sl renormalized 
Chern-Gauss-Bonnet formulas} by relating the total integral of such 
characteristic 
polynomials with the expected characteristic numbers and some
invariants of the CR-structure at infinity.

For instance, in (complex) dimension $m=2$, the integral of the 
characteristic polynomial 
\begin{equation}
3\, c_2 - (c_1)^2
\end{equation}
in the curvature tensor $R$ of
the complete Kähler-Einstein metric of a pseudo-convex domain
$\Omega$ in $\CM^2$ is shown to converge, and it is the only one to 
behave this way. Moreover, one can prove that
\begin{equation}\label{eq:burns-epstein}
\int_{\Omega}  ( \, c_2 - \frac{1}{3}\,(c_1)^2 )(R) = \chi(\Omega) + \mu(\partial\Omega)
\end{equation}
where $\chi$ is the Euler characteristic and $\mu$ is the 
{\sl Burns-Epstein invariant} of the CR-structure
of $\partial\Omega$ \cite{BurEps88}.

In \cite{Biq00}, the first author showed that a lot
of not necessarily integrable pseudoconvex
CR-structures on odd-dimensional spheres 
$\cerc^{2m-1}$ (including a neighborhood of the standard structure) 
may be filled in by complete Einstein metrics on the ball $\BM^{2m}$. 
Manifolds of (real) dimension $3$ ({\sl i.e.} $m=2$, as above) have the
special feature that there is no integrability condition for CR-structures.
However, it is well-known that a lot of such structures, even close to
the standard one, cannot be
obtained as the boundary of a complex domain, and hence, even in
this dimension the filling Einstein metrics cannot
be reduced to the classical complete Kähler-Einstein metrics of
pseudo-convex domains \cite{CheYau80}.

Nevertheless, the integrability of the CR structure at infinity supports
the idea that, in (real) dimension $2m=4$, the {\sl asymptotically complex 
hyperbolic} Einstein metrics (ACHE, in short) of \cite{Biq00} should retain 
some of the features of the Kähler situation. In this paper, we show
that this is indeed 
the case as far as renormalized characteristic classes are
concerned. More precisely, for any \ache manifold $(M,g)$ of (real) 
dimension~$4$, the special
combination of the norms of various parts of its curvature
\begin{equation}\label{eq:combiW}
\frac{1}{8\pi^2}\,
\left( 3\, |W^-|^2 - |W^+|^2 + \frac{1}{24}\,\Scal^2 \right)
\end{equation}
has convergent integral. Of course, this 
is the integrand of the characteristic 
class $\chi - 3\tau$ (where $\chi$ is the Euler characteristic and
$\tau$ the signature) which, in Kähler-Einstein geometry, may be 
rewritten as
$3c_2 - (c_1)^2 $.  Thus, our main result reads:

\medskip

\begin{maintheo}\label{maintheo}
Let $(M,g)$ be an asymptotically complex hyperbolic, Einstein, manifold
of dimension $4$. Then 
\begin{equation}\label{eq:formulefinale}
\frac{1}{8\pi^2}\,\int_M \left( 3\, |W^-|^2 - |W^+|^2 + 
\frac{1}{24}\,\Scal^2 \right)
\end{equation}
converges, and provides an invariant of the asymptotically complex
hyperbolic Einstein metric, which we call the Burns-Epstein invariant
of $g$.
\end{maintheo}

\medskip

Furthermore, one can hope to relate its values to the characteristic 
numbers and to some invariants of the CR structure of the boundary at 
infinity $\pb M$. We define such an invariant, proving:

\medskip

\begin{btheo} \label{btheo}
There is a global invariant $\nu(X)$ defined for any compact strictly
pseudoconvex 3-dimensional CR manifold $X$, such that if $X=\pb M$
is the CR structure at infinity induced by an asymptotically complex
hyperbolic, Einstein, metric on $M$, then
\begin{equation}\label{eq:formule2}
\frac{1}{8\pi^2}\,\int_M \left( 3\, |W^-|^2 - |W^+|^2 + 
\frac{1}{24}\,\Scal^2 \right) = \chi(M) - 3\,\tau(M) + \nu(X).
\end{equation}
\end{btheo}

\medskip

We now relate our invariant $\nu$ to the Burns-Epstein invariant
$\mu$. Remind that $\mu(J)$ is defined only in the case where the CR
structure $J$ on $X$ has trivial holomorphic bundle.
In the case where $M$ is Kähler-Einstein, and $\pb M$ has trivial
holomorphic bundle, one has the Burns-Epstein formula \cite{BurEps90b}
$$
\int_{M}c_{2}-\frac{1}{3}c_{1}^{2}=\chi(M)-\frac{1}{3}\bar{c}_{1}(M)^{2}+\mu(\pb M),
$$
where $\bar{c}_{1}(M)$ is a lifting of $c_{1}(M)$ in $H^{2}(M,\pb M)$.
The LHS is precisely one third the LHS in formula \eqref{eq:formule2},
but there is no clear relation between the topological terms on the
right, so we can only deduce that in that case $\nu(J)-3\mu(J)$ is
given by some topological term depending of the filling.
Of course, for the boundary $X$ of a domain in $\CM^2$, we get
immediately $\nu(X)=3\mu(X)+2$.

This problem can be avoided by considering the relative version $\mu(J,J')$
introduced by Cheng and Lee \cite{CL90}, defined now
for any CR structures $J$ and $J'$, and such that $\mu(J,J')=\mu(J)-\mu(J')$
when $J$ and $J'$ have trivial holomorphic bundle.
Then one can prove the following result.

\begin{ctheo}\label{ctheo}
For any CR structures $J$ and $J'$ on $X^{3}$ with the same underlying
contact structure, one has $$ \nu(J)-\nu(J')=3\mu(J,J'). $$
\end{ctheo}

This is a result which depends only the variations on the invariants
$\mu$ and $\nu$, and these variations are given by integration of
local terms.

By contrast, the invariant $\nu$ itself is defined as a
``renormalization'' of the $\eta$ invariants of $X$ for a sequence of
metrics converging to the Carnot-Carathéodory metric defined by the CR
structure, hence involving non local terms.
This explains why it is difficult to relate it to the $\mu$ invariant
of Burns-Epstein when the latter is defined.

It would be interesting to have a direct definition of $\nu$ on the CR
manifold, analogous to the spectral definition of the $\eta$ invariant,
instead of the definition by this limiting process.
Once this issue has been settled, our theorem \ref{btheo} will stand
as an analogue of 
the signature formula proven by N.~Hitchin for asymptotically {\sl real} 
hyperbolic Einstein metrics in dimension~$4$~\cite{Hit97}. A not-so-close 
analogue is the Gauss-Bonnet formula discovered by M.~T.~Anderson 
in the real case, which includes a contribution of the so-called
{\sl renormalized volume}, a non-local interior contribution~\cite{And01}. 

\medskip

The situation is by far less clear in higher dimension, as there is less
proximity between the general asymptotically complex hyperbolic case
and the very special Kähler-Einstein situation. If renormalized
characteristic classes were to exist, it seems difficult to extend the
methods of proof used in this paper to that case, as they use heavily 
the integrability of the CR-structure at infinity and the connection
this implies with complex geometry.

\bigskip

The paper is organized as follows.
In section~\ref{sec:asymptotics}, we study the asymptotics of an
asymptotically hyperbolic metric, various adapted connections and
their curvatures. This is used to refine our model at infinity, using 
Kähler geometry: in section~\ref{sec:approx}, given a CR-structure 
on any manifold $X^3$, we construct an explicit approximate metric,
which is Kähler-Einstein up to a high order. 
This involves generalizing the classical Fefferman procedure
\cite{Fef76} for complex domains to the abstract CR setting.
In section~\ref{sec:gauge}, we compare an arbitrary \ache metric $g$ 
with the approximate Kähler-Einstein metric $\bg$ built in the 
previous sections from the same CR-structure at infinity. Up to gauge 
modification (action of the diffeomorphism group), the former is shown 
to be a good approximation of the latter, up to a precise order; this 
is done in section~\ref{sec:da}. A useful output is the
explicit derivation of all the formally determined terms in the
asymptotic expansion of an \ache metric in dimension~$4$.

Unfortunately, this is not good enough to show that the characteristic 
polynomial (\ref{eq:combiW}) in the curvature of $g$ has convergent 
integral since the highest-order  term in the difference between $g$ and 
$\bg$ might cause divergence. In section~\ref{sec:weyl}, we show 
that the integrals converge by a direct method. The Einstein condition
implies that both half-Weyl tensors are harmonic. Thus, its negative part 
has fast enough decay to imply convergence of the $|W^-|^2$-term, whereas
the positive part can be compared to the positive part of the 
neighboring approximately Kähler-Einstein metric $\bg$. As 
$|W^+|^2 - \frac{1}{24}\,\Scal^2$ vanishes at least to high order for such 
a metric, the integral of the same term for $g$ can be shown to converge. 

In the following section~\ref{sec:comput}, we attack the task of computing 
the value of the integral. We transform it into a boundary integral by 
using the formulas for characteristic classes of manifolds with boundary 
and we consider the effect of 
the highest order term in $g-\bg$. Although it might contribute in the 
limit at infinity, a careful computation shows it is not the case. Our 
main result is then that the boundary term may be computed by using the
formally determined terms of the Kähler-Einstein metric $\bg$ rather than 
the \ache metric $g$. This enables us to define the invariant $\nu(\pb M)$
at infinity in full generality, for any CR structure on a $3$-dimensional
manifold, and to prove theorem \ref{btheo}.
Finally, in section \ref{sec:nu-variation}, we investigate the
relations between our invariant and the Burns-Epstein invariant $\mu$.
We compute the variation of $\nu$ with respect to the complex structure
and compare it to that of $\mu$, leading to the proof of theorem \ref{ctheo}.

\bigskip

\subsection*{Notations} We shall consider hereafter noncompact 
Riemannian manifolds, usually denoted as $(M,g)$. 
Their covariant (Levi-Civita) derivatives 
on any tensor bundle will be denoted by the symbol $\nabla$ and the 
Riemann, Ricci and scalar curvatures by $R^g$, $\Ric^g$ and $\Scal^g$ 
(the superscript being sometimes dropped if there is no possible 
confusion).
The sign convention on curvature is 
$ R_{X,Y} = [\nabla_X,\nabla_Y] - \nabla_{[X,Y]} $.  
The divergence
$\nabla^* = - \tr \nabla$ is the adjoint of $\nabla$, and the (rough) 
Laplacian operator on functions or tensors is then $\Delta^g = \nabla^* 
\nabla$.

\bigskip

\section{Asymptotics of ACH metrics}
\label{sec:asymptotics}

On the 3-sphere \( S^{3}\subset \mathbb{R}^{4} \), we denote
by \( \eta_0  \) the standard contact form, and \( \gamma_0  \) the
metric induced on the contact distribution \( \ker \eta  \).
The complex hyperbolic metric on \( \chii  \), with holomorphic
sectional curvature normalized to \( -1 \), is given in polar coordinates
by \begin{equation}
\label{for-ch2}
g_{\chii }=dr^{2}+4\sinh ^{2}(r)\eta ^{2}+4\sinh ^{2}(\frac{r}{2})\gamma .
\end{equation}

More generally, given any pseudo-convex
CR-structure on a 3-manifold \( X^{3} \),
a choice of compatible contact form \( \eta  \) induces a choice of 
metric \( \gamma (\cdot ,\cdot )=d\eta (\cdot ,J\cdot ) \)~; from this
one can build an {\sl asymptotically complex hyperbolic} 
metric on a neighborhood \( M=[R,\infty )× X \) of \(X\), 
\begin{equation}
\label{for-gACH}
g_0=dr^{2}+\exp^{2r}\eta ^{2}+\exp^{r}\gamma .
\end{equation}
As explained in \cite[I.1.B]{Biq00}, the curvature of \( g_0 \) is 
approximated by the curvature of \( g_{\chii } \) up to order
\( O(\exp^{-r/2}) \). This motivates the terminology {\sl asymptotically
complex hyperbolic metric} (ACH in short); a more general and precise 
definition will be given at the end of the current section. 
Note here the order \( O(\exp^{-r/2}) \) instead
of \( O(\exp^{-r}) \) in \cite{Biq00}, because we have normalized the
holomorphic sectional curvature to \( -1 \) instead of \( -4 \).

We will need later some calculations on the asymptotics of the metric
\( g_0 \). We let \( R \) be the Reeb vector field of the contact form 
\(\eta\), defined by 
\[ R\lrcorner \eta = 1, \ \ R\lrcorner d\eta =0 ,\]
and consider some unit vector field \( h \) in the contact distribution
\( H=\ker \eta  \). The CR-structure yields an almost complex structure
\( J \) on \( H \), which can be extended to an almost complex structure
on \( M \) by taking \( J\partial _{r}=\exp^{-r}R \). We may then consider 
an adapted \( g_0 \)-orthonormal frame 
\( (\partial _{r},\exp^{-r}R,\exp^{-r/2}h,\exp^{-r/2}Jh) \).
Finally, we denote by \( \nabla ^{W} \) the Webster connection on 
\( X \) determined by the choice of the contact form $\eta$.
Its torsion induces a trace-free symmetric
endomorphism of \( H \), 
\begin{equation}
\label{for-torW}
T_{R,\cdot }=\left( \begin{array}{cc}
\alpha  & \beta \\
\beta  & -\alpha 
\end{array}\right) .
\end{equation}
We can extend the Webster connection 
to \( M \) as a \( g_0 \)-unitary connection 
\( \widetilde{\nabla }^{W} \)by defining 
\[
\widetilde{\nabla }^{W}\partial _{r}=\widetilde{\nabla }^{W}\exp^{-r}R=0,
\qquad \widetilde{\nabla }^{W}_{\partial _{r}}\exp^{-\frac{r}{2}}h=0
\textrm{ for }h\in H.
\]
\begin{lem}\label{lem-nabla}
The Levi-Civita connection of \( g_0 \) is 
\( \nabla =\widetilde{\nabla }^{W}+a \),
where \( a \) is a 1-form with values in the endomorphisms of \( TM \)
defined in the \( g_0\)-adapted frame 
\( (\partial _{r},\exp^{-r}R,\exp^{-r/2}h,\exp^{-r/2}Jh) \)
by \begin{eqnarray*}
a_{\partial _{r}} & = & 0,\quad a_{\exp^{-r}R}=\left( \begin{array}{cccc}
0 & -1 &  & \\
1 & 0 &  & \\
 &  & 0 & -\frac{1}{2}\\
 &  & \frac{1}{2} & 0
\end{array}\right) ,\\
a_{\exp^{-\frac{r}{2}}h} & = & \left( \begin{array}{cccc}
 &  & -\frac{1}{2} & 0\\
 &  & -\exp^{-r}\alpha  & -\frac{1}{2}-\exp^{-r}\beta \\
\frac{1}{2} & \exp^{-r}\alpha  &  & \\
0 & \frac{1}{2}+\exp^{-r}\beta  &  & 
\end{array}\right) ,\\
a_{\exp^{-\frac{r}{2}}Jh} & = & \left( \begin{array}{cccc}
 &  & 0 & -\frac{1}{2}\\
 &  & \frac{1}{2}-\exp^{-r}\beta  & \exp^{-r}\alpha \\
0 & -\frac{1}{2}+\exp^{-r}\beta  &  & \\
\frac{1}{2} & -\exp^{-r}\alpha  &  & 
\end{array}\right) .
\end{eqnarray*}
\end{lem}

\medskip

\begin{rem}\label{rem:reexpressed}
The result of the previous lemma may be re-expressed in the following way: 
\begin{equation}\label{for-a01}
a=a_{0}+\exp^{-r}a_{1}
\end{equation}
where \( a_{0} \) is a $1$-form having the same coefficients as it has
in \( \chii  \) (in particular,
\( a_{0} \) commutes with \( J \)); the correction term \( a_{1} \)
depends on the torsion of the Webster connection. This stands also true
for all derivatives.
\end{rem}

\medskip

\begin{proof}
The proof is a straightforward calculation, using the fact that, for
an orthonormal frame \( (\xi _{i}) \), the Levi-Civita connection
can be computed by the formula \[
\left\langle \xi _{k},\nabla _{\xi _{j}}\xi _{i}\right\rangle 
=-\frac{1}{2}
\left( \left\langle [\xi _{i},\xi _{k}],\xi _{j}\right\rangle 
+\left\langle [\xi _{j},\xi _{k}],\xi _{i}\right\rangle 
+\left\langle [\xi _{i},\xi _{j}],\xi _{k}\right\rangle \right) .
\]
All the brackets can be expressed in terms of the Webster connection
and its torsion, since for \( h,h'\in H \), one has 
\[
T_{h,h'}=d\eta (h,h')R,
\]
and the other components of the torsion are given by (\ref{for-torW}).
\end{proof}

\medskip

It is proven in \cite{Biq00} that the curvature of \( g_0 \) is approximated 
up to order \( O(\exp^{-r/2}) \) by the curvature of the model space 
\( \chii  \); therefore, \( \Ric^{g_0}=-\frac{3}{2}g_0+O(\exp^{-r/2}) \),
hence $g_0$ is a first approximation solution to the Einstein equation.
More is actually true: the correction
term in lemma \ref{lem-nabla} is better than expected and decays like 
\( O(\exp^{-r}) \) instead of \( O(\exp^{-r/2}) \). In order to see this,
we first deduce from the previous Lemma:

\begin{cor}
\label{cor-nablaa}
The Levi-Civita connection of \( g_0 \) and the
Webster connection are related by \( \nabla =\widetilde{\nabla }^{W}+a, \)
with 
\[
\widetilde{\nabla }^{W}a=O(\exp^{-r}),
\]
where \( O(\exp^{-r}) \) is taken with respect to \( g_0 \). The same is true
for all derivatives: they all are of order \( O(\exp^{-r}) \).
\end{cor}

\begin{proof}
From the decomposition (\ref{for-a01}) of \( a \), the only thing
to prove is \( \widetilde{\nabla }^{W}a_{0}=0 \). Actually one can give
more intrinsic formulas for \( a_{0} \). Denote by \( \widetilde{g} \)
the hermitian product on \( TM \) induced by \( g_0 \) and \( J \),
which is \( \mathbb{C} \)-antilinear in its first variable. Then
\[
(a_{0})_{\exp^{-r}R}|_{\left\langle \partial _{r},\exp^{-r}R\right\rangle }
=J,\quad (a_{0})_{\exp^{-r}R}|_{H}=\frac{1}{2}J,
\]
and for any \( h\in H \), 
\[
(a_{0})_{\exp^{-r/2}h}\xi =\frac{1}{2}
\left( \widetilde{g}(\partial _{r},\xi )\exp^{-r/2}h-\widetilde{g}(\exp^{-r/2}h,\xi )
\partial _{r}\right) .
\]
Since \( J \) is parallel for \( \widetilde{\nabla }^{W} \), it follows
that \( a_{0} \) is parallel too, and this ends the proof.
\end{proof}

\begin{cor}\label{cor-Ricg}
The curvature of the metric \( g_0 \) defined by (\ref{for-gACH})
is approximated up to order \( O(\exp^{-r}) \) by that of the complex
hyperbolic space ; in particular, one has 
\[
\Ric^{g}=-\frac{3}{2}g+O(\exp^{-r}).
\] 
\end{cor}

\medskip

\begin{rem}\label{remann}
Actually one can compute explicitely the term of order \( \exp^{-r} \):
it depends on the Webster scalar curvature and on the torsion. From this 
one may 
write down a first correction (at order \( \exp^{-r} \)) to \( g \) in 
order to get an approximate Einstein metric up to a better order.
In sections \ref{sec:approx} and \ref{sec:da}, we will prove an even better 
asymptotic expansion for a more precise choice of approximate Einstein metric.
\end{rem}

\medskip

\begin{proof}
This is also a consequence of Lemma \ref{lem-nabla}.
Indeed, we can write the curvature of the Levi-Civita connection as
\[
F=F(\widetilde{\nabla }^{W})+d_{\widetilde{\nabla }^{W}}a+\frac{1}{2}[a,a].
\]
Notice first that \( F(\widetilde{\nabla }^{W}) \) actually reduces
to \( F(\nabla ^{W}) \): since this is a smooth horizontal 2-form
on the boundary, it means that, with respect to the metric \( g_0 \),
we have 
\[
F(\widetilde{\nabla }^{W})=O(\exp^{-r}).
\]
Secondly, using (\ref{for-a01}) and corollary \ref{cor-nablaa},
we get
\[
d_{\widetilde{\nabla }^{W}}a+\frac{1}{2}[a,a]=a_0(T^W_{|H}) + 
\frac{1}{2}[a_{0},a_{0}]+O(\exp^{-r})\]
where $T^W$ is the torsion of the Webster connection.
Therefore, we conclude that \[
F=a_0(T^W_{|H}) + \frac{1}{2}[a_{0},a_{0}]+O(\exp^{-r}).
\]
The form \( a_{0} \) has constant coefficients equal to those of the model: 
the term \( a_0(T^W_{|H}) + \frac{1}{2}[a_{0},a_{0}] \) actually represents 
the curvature
of \( \chii  \), hence this formula implies that the curvature of \( g_0 \)
is approximated by that of complex hyperbolic space up to
order \( O(\exp^{-r}) \).
\end{proof}

\subsection*{Laplacians}
We close this section by another consequence of corollary
\ref{cor-nablaa}. Given any tensorial bundle \( E \) on \( X \),
we have on it transverse-Webster operators as $\nabla^W_R$, $\nabla^W_h$
(for $h$ section of $H$).
For handling transverse regularity questions, it is important to
understand the commutation of their extensions $\wn_R$ and $\wn_h$ with
the standard Laplacian \( \Delta_{g_0}  \) for the metric \( g_0 \) on \( M \).

Since \( M \) is asymptotically complex hyperbolic, 
the unit balls in \( M \) at infinity look
like the unit balls in \( \chii  \), and this enables us to define
in a standard way Hölder spaces \( C^{k,\alpha } \) for the metric
\( g \). We also need the weighted versions \[
C^{k,\alpha }_{\delta }=\exp^{-\delta r}C^{k,\alpha }.\]
 Clearly, the Laplacian for \( g \) defines an operator \[
\Delta:C^{k+2,\alpha }_{\delta }\to C^{k,\alpha }_{\delta },\]
and the tangential derivatives give operators
\[
\wn_h : C^{k+1,\alpha }_{\delta }\to C^{k,\alpha }_{\delta -\frac{1}{2}}, \ 
\wn_R : C^{k+1,\alpha }_{\delta }\to C^{k,\alpha }_{\delta -1}, 
\]
since the norm of \( h \) for
the metric \( g \) is \( \exp^{r/2} \), resp $\exp^r$ for $R$. Hence the
bracket \( [\Delta , \mathcal{L} ] \) loses {\it a priori} some weights
if $\mathcal{L}$ is any of the above transverse-Webster operators. 
Although this might put us in a very bad shape, we shall see now that
the case is actually better.

\medskip

We denote by $\cQ_1$ the algebra of differential operators generated
by $\wn_{\partial_r}$, $\wn_{\exp^{-r}R}$, $\wn_{\exp^{-\frac{r}{2}}h}$ 
($h$ smooth section of $H$ on $X$),
and linear (zero order) operators $f$ on $M$ such that,
for any integer $k$ and sections $h_1,\dots,h_k$ of $H$, 
$\wn_{h_1}\dots\wn_{h_k} f$ is bounded and, for any $\ell>0$, 
$\wn_{h_1}\dots\wn_{h_k}(\wn_{\partial_r})^{\ell}f = O(\exp^{-\frac{r}{2}})$
(in all what follows, these will be series in $\exp^{-\frac{r}{2}}$
whose coefficients are smooth functions on the boundary
at infinity $X$). 
Let $\cQ_0$ be the subspace of $\cQ_1$ containing operators with only
tangential derivatives and no zero order term.
Clearly,
\[ \cQ \ = \ \cQ_0 + \exp^{-\frac{r}{2}}\,\cQ_1 \]
is an algebra of differential operators.
The main interest of this algebra lies in its commutation properties with 
the transverse-Webster operators or the metric Laplacian of $g_0$:

\begin{lem}\label{lem:commut} One has, for $h$ unit section of $H$ on $X$, 
\begin{equation}\label{eq:commut}\begin{split} 
& [\wn_{\partial_r},\wn_h ] = 0, \ [\wn_{\partial_r},\wn_R ] = 0, \\
& [\wn_h,\wn_{\exp^{-r}R} ] \in \exp^{-\frac{r}{2}}\cQ_1, \\
& [\wn_h,\wn_{\exp^{-\frac{r}{2}}Jh}] = - \wn_{ \exp^{-\frac{r}{2}}R} \
\textrm{ mod } \cQ .  
\end{split}\end{equation}
As a consequence, 
\begin{equation}\label{eq:commut2}
[\Delta_{g_0}, \wn_R] \in \cQ \ \textrm{ and } 
\ [\Delta_{g_0}, \wn_h] = - 2 \,\nabla_{\exp^{-\frac{r}{2}}Jh}
\wn_{ \exp^{-\frac{r}{2}}R} \ \textrm{ mod } \cQ .\end{equation}
\end{lem}

\begin{proof}
The first commutations properties can be easily checked from the form
of the Webster 
curvature and contact properties of the distribution  $H$. It remains to
prove the (more interesting) commutation with the metric Laplacian. In 
a \( g_0 \)-orthonormal frame $ (e_{i}) = (\partial_r,\exp^{-r}R,
\exp^{-\frac{r}{2}}h,\exp^{-\frac{r}{2}}Jh)$, the Laplacian is 
\begin{equation}\label{for-Delta}
\Delta \  = \ -\sum \left( (\nabla _{e_{i}})^{2}
-\nabla _{\nabla _{e_{i}}e_{i}} \right)
\end{equation}
where the last line is obtained from Lemma \ref{lem-nabla}: $\nabla = \wn + a_0
+ \exp^{-r}a_1$ with $a_1\in \cQ_1$ and $\wn a_0 =0$.
From the previous commutations, one sees that the only term outside 
$\cQ$ in $[\wn_h,\Delta_{g_0}]$ is
\begin{equation*}\begin{split}
[\wn_h,-\left(\nabla_{\exp^{-\frac{r}{2}}Jh}\right)^2] & = 
- \nabla_{\exp^{-\frac{r}{2}}Jh}[\wn_h,\nabla_{\exp^{-\frac{r}{2}}Jh}]
- [\wn_h,\nabla_{\exp^{-\frac{r}{2}}Jh}]\nabla_{\exp^{-\frac{r}{2}}Jh}\\
& = \wn_{\exp^{-\frac{r}{2}}R}\nabla_{\exp^{-\frac{r}{2}}Jh}
+ \nabla_{\exp^{-\frac{r}{2}}Jh}\wn_{\exp^{-\frac{r}{2}}R} \ \textrm{ mod }
\cQ \\
& = 2\,\nabla_{\exp^{-\frac{r}{2}}Jh}\wn_{\exp^{-\frac{r}{2}}R} + 
[\nabla_{\exp^{-\frac{r}{2}}R},\wn_{\exp^{-\frac{r}{2}}Jh}] \ \textrm{ mod }
\cQ \\
& = 2\,\nabla_{\exp^{-\frac{r}{2}}Jh}\wn_{\exp^{-\frac{r}{2}}R} \ \textrm{ mod }
\cQ\ .
\end{split}\end{equation*}
The commutations with $\wn_R$ are similar, but easier.
\end{proof}

\medskip

\subsection*{General definition of ACH(E) metrics} 
We now close this introductory section with a general definition of the
Riemannian manifolds which form the main objects of study of this paper.
Let $(X^3,\eta,J)$ 
be a pseudo-convex CR manifold, with 
associated metric $\gamma$ on the contact distribution. Any metric 
$g$ on a $4$-manifold $M$ such that the complement of a compact set $M-K$ 
is diffeomorphic to $[R,+\infty[× X$ and such that
\[ g - (dr^2 + \exp^{2r}\eta^2 + \exp^r\gamma)  \ \in \ C^{\infty}_{\delta}
\]
for some $\delta >0$ 
will be called an {\sl asymptotically complex hyperbolic}
(ACH in short) manifold. 

Moreover, $(M,g)$ is said to be ACHE if $g$ is an Einstein metric. To 
avoid any confusion, we insist on the fact that ACHE metrics actually
are solutions of the Einstein equations; the word `asymptotic' in the 
definition only refers to the complex hyperbolic-like behavior at infinity:
each should look like complex hyperbolic space.

From now on and for sake of simplicity, we shall restrict ourselves to 
smooth ACH metrics induced by smooth CR structures at infinity, and 
weighted decay control on {\sl all} derivatives, as it is the case in the
definition above. This is because we know from the work of the 
first author \cite{Biq00} that such metrics can be obtained on the ball 
from any smooth CR structure close to the standard structure on $\cerc^3$. 
This provides us with a very large set of metrics to which our results can 
be applied. Their domain of validity may likely be pushed further to 
include finite differentiability assumptions only, although this might 
require a slightly more technical treatment.

\bigskip

\section{Approximately Kähler-Einstein metrics}\label{sec:approx}

\medskip

Let us recall that any strictly pseudoconvex domain in \( \mathbb{C}^{2} \)
bears in its interior a complete Kähler-Einstein metric: the Cheng-Yau
metric \cite{CheYau80}, whose asymptotics are similar to those of the
complex hyperbolic metric (\ref{for-ch2}). Fefferman \cite{Fef76}
has given a formal high order asymptotic expansion for such a metric, and 
Lee-Melrose \cite{LeeMel82} proved the complete asymptotic expansion.

More generally, given any CR-structure on a 3-manifold \( X^{3} \),
we have at hand the asymptotically complex hyperbolic metric 
(\ref{for-gACH}) of the previous section,
given on a neighborhood \( [R,\infty )× X \) by 
\begin{equation}
\label{for-g}
g=dr^{2}+\exp^{2r}\eta ^{2}+\exp^{r}\gamma .
\end{equation}
We will now modify this metric in order to get a Kähler-Einstein
metric, at least up to a very high order. As in general \( X \) is not 
embedded in \( \mathbb{C}^{2} \), we cannot use directly
Fefferman's formal construction.

We begin by the construction of the complex structure.
This is well-known, but we need the calculation of the first terms.
\begin{prop}
\label{pro-formel}Given any CR manifold \( X^{3} \), one can construct
in a neighborhood of \( X \) a formal integrable complex structure \( J \).
\end{prop}

\begin{proof}
Denote by \( J_{0} \) the complex structure on the contact 
distribution, and \( R \) the Reeb vector field associated to the contact 
form \( \eta  \). Moreover, let \( H_{1,0} \) be the (1,0)-vectors 
in \( H \) and \( H^{0,1} \) be the (0,1)-forms. 

We consider now \( M=X× [R,\infty ) \), and we start from the initial 
almost complex structure on \( M \) defined by the formulas 
\[
J_{0}\partial _{r}=\exp^{-r}R,\quad J_{0}|_{H}=J_{X}.
\]
We now seek a complex structure \( J \) differing from \( J_{0} \) only
on \( H \): the difference is parameterized by a tensor 
\( \phi \in H^{0,1}\otimes H_{1,0} \),
such that 
\[
T^{J}_{0,1}=\{X+\phi _{X},\, X\in T^{J_{0}}_{0,1}\}.
\]
The integrability condition 
\( [T^{J}_{0,1},T^{J}_{0,1}]\subset T^{J}_{0,1} \)
becomes, in terms of \( \phi  \) and an arbitrary vector 
\( h\in H_{0,1} \),
\[
[\partial _{r}+i\exp^{-r}R,h+\phi _{h}]\in (1+\phi )T^{J_{0}}_{0,1}.
\]
We perform the calculation at a point \( x\in X \), where we suppose 
that \( \nabla^W h(x)=0 \) for the Webster connection $\nabla^W$, 
and therefore
 \( [R,h]=-T_{R,h} \), where \( T \) is the torsion of the Webster 
connection. As \( T_{R,\cdot } \) is anti-\( J_{0} \)-linear on \( H \)
and therefore defines a map \( H_{0,1}\to H_{1,0} \), we
get 
\[
[\partial _{r}+i\exp^{-r}R,h+\phi _{h}]=-i\exp^{-r}T_{R,h}
+\partial _{r}\phi _{h}+i\exp^{-r}[R,\phi _{h}] \ .
\]
Using \( [R,\phi _{h}]=\nabla^W _{R}\phi _{h}-T_{R,\phi _{h}} \),
the integrability condition can be transformed into 
\[
-i\exp^{-r}T_{R,h}+\partial _{r}\phi _{h}+i\exp^{-r}\nabla^W _{R}\phi _{h}
=-i\exp^{-r}\phi _{T_{R,\phi _{h}}},
\]
which we rewrite finally as
\begin{equation}
\label{eq-int}
-\partial _{r}\phi _{h}=i\exp^{-r}\left( -T_{R,h}+\nabla^W _{R}\phi _{h}
+\phi _{T_{R,\phi _{h}}}\right) .
\end{equation}

Now it becomes clear that we can solve (\ref{eq-int}) by a formal
series \( \phi =\sum _{j\geq 1}\phi _{j}\exp^{-jr} \). Indeed, suppose
we have a solution up to order \( k-1 \), then the r.h.s. of (\ref{eq-int}),
computed for \( \sum _{j<k}\phi _{j}\exp^{-jr}, \) is at least of
order \( k \), {\sl i.e.} of the form \( \sum _{j\geq k}\exp^{-jr}\psi _{j} \),
and we can solve the equation at order \( k \) by letting 
\[ \phi _{k}=\frac{1}{k}\psi _{k} \]
and this proves our claim.
\end{proof}

\medskip

\begin{rem}
Formula (\ref{eq-int}) enables us to give easily an explicit formula
for the first terms of $\phi$. For example, we have clearly $\phi_1
=-iT_{R,\cdot}$.
Reintroducing this into (\ref{eq-int}), we get that $2\phi_2=
i\nabla_R \phi_1=\nabla_R T_{R,\cdot}$. As a result,
\begin{equation}\label{dv-phi}
\phi = -iT_{R,\cdot} \exp^{-r} + \frac{1}{2} \nabla_R T_{R,\cdot} \exp^{-2r}
       +\cdots
\end{equation}
yields the beginning of the series for the complex structure.
\end{rem}

\medskip

We will now construct the approximate Kähler-Einstein metric. Before
stating the theorem, we need to recall some formalism for the Webster
connection on the CR manifold $X$.
We work in a local coframe $(\eta,\ti,\tib)$ such that
$$ d\eta=i\ti\land\tib. $$
The Webster connection form is a purely imaginary 1-form $\oii$, and
the Webster torsion $\taui=\taui_{\bar{1}}\tib$ is a (0,1)-form: 
they are defined by
$$ d\ti = \ti \land \oii + \eta \land \taui. $$
The Webster curvature $R$ is defined by
$$ d\oii = R\ti\land\tib
          + \tauib_{,\bar{1}}\land\eta - \taui_{,1}\land\eta\ , $$
where $D_{,1}$ (resp. $D_{,\bar 1}$, or $D_{,0}$) denotes covariant 
derivative 
(with respect to the Webster connection) of the tensor field $D$
in the direction of the $(1,0)$ vector dual to $\theta^1$ 
(resp. $\theta^{\bar 1}$, or the contact form $\eta$).
The canonical bundle of the CR structure (generated by
$\ti\land\eta$) is a natural CR holomorphic bundle; it has a
canonical connection with curvature $\Omega$ satisfying 
$\Omega\land\eta=0$.
More precisely,
$$ \Omega = \eta\land\left((iR_{,1}\ti+\tauib_{,\bar{1}})
           +(iR_{,\bar{1}}\tib-\taui_{,1})\right). $$
From this point of view, the complex structure of the filling complex
structure defined by (\ref{dv-phi}) is determined by the (1,0)-forms
\begin{equation}\label{eq:vtheta}\begin{split}
\vt&=\exp^{-r}dr+i \eta, \\
\vti&=\ti-\phi\lrcorner\ti
  =\ti+i\exp^{-r}\taui-\frac{1}{2}\exp^{-2r}\taui_{,0}+\cdots
\end{split}\end{equation}
Therefore, if we see $\Omega$ as a 2-form in the interior, the leading
term of its (1,1)-part is
$$
\Omega_0=\frac{i}{2}
\left(\vtb\land(iR_{,1}\ti+\tauib_{,\bar{1}})
     -\vt\land(iR_{,\bar{1}}\tib-\taui_{,1})
\right).
$$
This form decreases as $O(\exp^{-\frac{3}{2}r})$, but is not closed, since
an easy calculation gives
$$
d\Omega_0=-\frac{i}{2}
(\Delta R-i(\taui_{\bar{1},11}-\tauib_{1,\bar{1}\bar{1}}))
\exp^{-r}dr\land d\eta + O(\exp^{-\frac{5}{2}r}),
$$
where $\Delta$ is the Webster Laplacian.
This leads us to define a form
$$
\tilde{\Omega}=\Omega_0-\frac{i}{2}
(\Delta R-i(\taui_{\bar{1},11}-\tauib_{1,\bar{1}\bar{1}}))
\exp^{-r} d\eta ,
$$
which now satisfies
$$ d\tilde{\Omega}=O(\exp^{-\frac{5}{2}r}). $$

\medskip

\begin{theo}\label{theo:omegaKE}
In a neighborhood of $X$ there exists a formal Kähler ACH metric $\bg$, 
which is Einstein up to order $O(\exp^{-\frac{5}{2}r})$, that is
$$\Ric^{\bg}=-\frac{3}{2}\bg+O(\exp^{-\frac{5}{2}r}).$$
Moreover, the Kähler form $\omega$ of the metric is given explicitely, 
up to order $\frac{5}{2}$, by: 
$$
\omega=i\partial\dbar f+\frac{4}{3}i\tilde{\Omega}
       +O(\exp^{-\frac{5}{2}r}),
$$
and
$$
f=2r+R\exp^{-r}
 -\frac{2}{3}\left(\frac{R^2}{4}-|\tau|^2-\frac{\Delta R}{6}
 +\frac{2i}{3}(\taui_{\bar{1},11}-\tauib_{1,\bar{1}\bar{1}})\right)
  \exp^{-2r}.
$$
\end{theo}

\medskip

\begin{cor}\label{cor:omega-explicit} Up to order $\frac{5}{2}$, one obtains
explicitely  
\begin{equation*}\begin{split}
\omega \ = \ & \exp^r \left( dr\land\eta + d\eta\right) \ 
- \ \frac{R}{2}\, d\eta \\
&\ + \frac{4}{3}i\tilde{\Omega}
\ + \ \frac{i}{2}\left(R_{,1}\vtb\land\ti
\ - \ R_{,\bar 1}\vt\land\tib\right)
\ - \ \frac{\Delta R}{2} \exp^{-r}d\eta\\
&\ - \ \frac{2}{3} \left(\frac{R^2}{4}-|\tau|^2-\frac{\Delta R}{6}
 +\frac{2i}{3}(\taui_{\bar{1},11}-\tauib_{1,\bar{1}\bar{1}})\right)
 \exp^{-r}(dr\land\eta-d\eta),
\end{split}\end{equation*}
where $\vt$, $\vti$ are defined in formulas \upn{(}\ref{eq:vtheta}\upn{)}.
\end{cor}

\smallskip

\begin{rem}
Actually it is easy to prove that there is a formal development,
determined up to order 3, for an additional potential for a solution,
as in Fefferman's classical work \cite{Fef76}.
However, the development up to order $\frac{5}{2}$ obtained above will 
be enough for our needs.
\end{rem}

\smallskip

\begin{proof}
First we have
$ \dbar(2r) = dr - i \exp^r \eta $,
and therefore
\begin{equation}\label{dd2r}
i\partial\dbar(2r)=id\dbar(2r)=\exp^r(dr\land\eta+d\eta).
\end{equation}
Hence the leading term for $g$ actually gives the ACH metric
(\ref{for-g}).
Now continue the calculation:
$$
d(R\exp^{-r})=\exp^{-r}(-Rdr+R_{,0}\eta+R_{,1}\ti+R_{,\bar{1}}\tib).
$$
Remark that $\tib=\vtib+O(\exp^{-\frac{3}{2}r})$, where of course the
$O(\exp^{-\frac{3}{2}r})$ refers to the ACH metric (\ref{for-g}).
Therefore
$$ \dbar(R\exp^{-r})=
\frac{1}{2}(-R+i\exp^{-r}R_{,0})\vtb+\exp^{-r}R_{,\bar{1}}\tib
+O(\exp^{-\frac{5}{2}r}). $$
We deduce, keeping terms only up to order $\exp^{-2r}$,
\begin{multline*}
\partial\dbar(R\exp^{-r})=
-(R_{,0}\eta+R_{,\bar{1}}\tib+R_{,1}\ti)\frac{\vtb}{2}
+\frac{R}{2}id\eta
+\frac{R_{,0}}{2}\exp^{-r}(-dr\land\eta+d\eta)\\
-\exp^{-r}dr\land R_{,\bar{1}}\tib
+\exp^{-r}R_{,\bar{1}1}\ti\land\tib+O(\exp^{-\frac{5}{2}r})
\end{multline*}
and therefore
\begin{equation}\label{ddRer}\begin{split}
i\partial\dbar(R\exp^{-r})\ = \ & \frac{1}{2}\left(
-R d\eta - i\vt\land R_{,\bar{1}}\tib + i\vtb\land R_{,1}\ti
- \exp^{-r}(\Delta R) d\eta \right) \\ & \ + \ O(\exp^{-\frac{5}{2}r}).
\end{split}
\end{equation}
We consider the metric Kähler metric $g_0$ with Kähler form
$$ \omega_0 = i\partial\dbar(2r+R\exp^{-r}). $$
Let us calculate its Ricci tensor. We have a section
$$ \sigma = \vt \land \vti $$
of the canonical bundle.
Remark that
$$ i\vti\land\vtib
 = i(\ti\land\tib+\exp^{-2r}\taui\land\tauib)+O(\exp^{-4r}).
$$
It follows that, with respect to $g_0$, we have
$$|\sigma|^2 = \exp^{-2r}
\frac{1-\exp^{-2r}|\taui|^2}{e^r-\frac{1}{2}(R+\Delta R \exp^{-r})}
(1+O(\exp^{-\frac{5}{2}r})),
$$
and
\begin{equation}\label{lnsigma2}
\ln|\sigma|^2=-3r+\frac{R}{2}\exp^{-r}+\Phi\exp^{-2r}
  +O(\exp^{-\frac{5}{2}r}),\quad 
\Phi=\frac{R^2}{4}-|\taui|^2+\frac{\Delta R}{2}.
\end{equation}
On the other hand, we have
$$ d\sigma=id\eta\land\vti
 -\vt\land(d\ti+d(i\exp^{-r}\taui-\frac{1}{2}\exp^{-2r}\taui_{,0})) ;
$$
the first term is zero, and the second becomes, keeping only terms of
type (2,1),
$$
d\sigma=-\vt\land(\ti\land\oii+\eta\land\taui+i\exp^{-r}d\taui)
+O(\exp^{-\frac{5}{2}r});
$$
replacing $\ti$ by $\vti-i\exp^{-r}\taui+O(\exp^{-\frac{5}{2}r})$, we get
$$
d\sigma=(-\oii-i\exp^{-r}\taui_{,1}+O(\exp^{-\frac{5}{2}r}))\land\sigma,
$$
giving us a ``connection'' form
\begin{equation}\label{omK}
\omega_K=-\oii-i\exp^{-r}(\taui_{,1}+\tauib_{,\bar{1}})
  +O(\exp^{-\frac{5}{2}r}).
\end{equation}
One deduces easily:
\begin{equation}\label{domK}
d\omega_K\ = \ -R\ti\land\tib+i\vt\land\taui_{,1}+i\vtb\land
\tauib_{,\bar{1}}
-i\exp^{-r}d(\taui_{,1}+\tauib_{,\bar{1}})+O(\exp^{-\frac{5}{2}r}).
\end{equation}
Putting together (\ref{dd2r}), (\ref{ddRer}), (\ref{lnsigma2}) and
(\ref{domK}), we calculate the Ricci form of the Kähler metric $g_0$:
\begin{align*}
\rho^{\omega_0} &= -i\dbar\partial\ln|\sigma|^2 -i d\omega_K \\
&= -\frac{3}{2}\exp^r (dr\land\eta+d\eta)
   +\frac{3}{4}R d\eta \\
&\quad +\frac{1}{4} \left(-i\vt\land R_{,\bar{1}}\tib + i\vtb\land R_{,1}\ti
- \exp^{-r}(\Delta R) d\eta \right) \\
&\quad +\Phi \exp^{-r}(dr\land\eta-d\eta) \\
&\quad +\vt\land\taui_{,1}+\vtb\land\tauib_{,\bar{1}}
-\exp^{-r}d(\taui_{,1}+\tauib_{,\bar{1}})+O(\exp^{-\frac{5}{2}r})
\end{align*}
Therefore
\begin{align*}
\rho^{\omega_0}+\frac{3}{2}\omega_0&=
-i\vt\land R_{,\bar{1}}\tib + i\vtb\land R_{,1}\ti
- \exp^{-r}\,(\Delta R) d\eta +\Phi \exp^{-r}(dr\land\eta-d\eta) \\
&\quad
+\vt\land\taui_{,1}+\vtb\land\tauib_{,\bar{1}}
-\exp^{-r}d(\taui_{,1}+\tauib_{,\bar{1}})+O(\exp^{-\frac{5}{2}r}) \\
&=-\vt\land(iR_{,\bar{1}}\tib-\taui_{,1})
  +\vtb\land(iR_{,1}\ti+\tauib_{,\bar{1}}) \\
&\quad +\Phi \exp^{-r}(dr\land\eta-d\eta)
       -\exp^{-r}((\Delta R)d\eta+d(\taui_{,1}+\tauib_{,\bar{1}}))+O(\exp^{-\frac{5}{2}r}) \\
&= -2i \tilde{\Omega} +\Phi \exp^{-r}(dr\land\eta-d\eta)
       +O(\exp^{-\frac{5}{2}r}).
\end{align*}
From this formula, we see that $g_0$ is a Kähler metric, which is
Einstein up to order $O(\exp^{-\frac{3}{2}r})$. The term of order
$\frac{3}{2}$ does not come from a potential in general, since
$\Omega$ represents the first Chern class of the canonical bundle of
the boundary.
However, we can modify easily $\omega_0$ in order to kill this
term. Indeed, observe that the order $\frac{3}{2}$ term in
$\tilde{\Omega}$, that is $\Omega_0$, is orthogonal to the
leading term (\ref{dd2r}) of the Kähler form (and also to the term of
order 1 in (\ref{ddRer})). This implies that if we define a new Kähler
form
$$\omega_1=\omega_0+\frac{4}{3}i\tilde{\Omega},$$
the metric on the canonical bundle is changed only at order
$O(\exp^{-2r})$ by the terms of order 2 in $\tilde{\Omega}$, hence
Ricci is modified only at the same order.
More precisely, we must add a term
$$ -\frac{2}{3}\big(\Delta R-i(\taui_{\bar{1},11}-\tauib_{1,\bar{1}\bar{1}})
               \big) $$
in formula \eqref{lnsigma2} for $\ln|\sigma|^2$, which amounts to
replace $\Phi$ by
\begin{align*}
\Phi_1&=
\Phi-\frac{2}{3}\big(\Delta R-i(\taui_{\bar{1},11}-\tauib_{1,\bar{1}\bar{1}})
                \big) \exp^{-r}\\
&=\frac{R^2}{4}-|\taui|^2-\frac{\Delta R}{6}+\frac{2i}{3}
(\taui_{\bar{1},11}-\tauib_{1,\bar{1}\bar{1}}).
\end{align*}
Putting things together, we now have
$$
\rho^{\omega_1}+\frac{3}{2}\omega_1=
\Phi_1 \exp^{-r}(dr\land\eta-d\eta)+O(\exp^{-\frac{5}{2}r}).
$$

It remains to kill the terms of order $2$. For this, observe that 
again the order~2 term $\Phi_1\exp^{-r}(dr\land\eta-d\eta)$ is
orthogonal to the leading term (\ref{dd2r}) of the Kähler form,
meaning that if we take
\begin{align*}
\omega_2&=\omega_{1}-\frac{2}{3}i\partial\dbar(\Phi_{1}\exp^{-2r})\\
  &=\omega_1-\frac{2}{3}\Phi_1\exp^{-r}(dr\land\eta-d\eta)+O(\exp^{-\frac{5}{2}r}),
\end{align*}
then Ricci is unchanged at order 2, so we obtain
$$ \rho^{\omega_2}+\frac{3}{2}\omega_2 = O(\exp^{-\frac{5}{2}r}). $$

Summarize all the corrections we have done (modulo
$O(\exp^{-\frac{5}{2}r})$) by the formula
$$
\omega_{2}=\omega_{0}-\frac{2}{3}
\left(\rho^{\omega_0}+\frac{3}{2}\omega_0-i\partial \dbar((\Phi+\Phi_{1})\exp^{-2r})
\right).
$$
This defines $\omega_2$ as a closed (1,1)-form, which is the Kähler form we
were looking for.
\end{proof}

\bigskip

\section{Putting the Einstein metric in an adequate gauge}\label{sec:gauge}

\bigskip

The goal of this section is to prove that, starting with an \ache metric
$g$ on a manifold $M^4$, with CR-structure $(\eta,H,J_0)$ at infinity 
$\pb M=X^3$, one can approach it by an approximate solution of the
Kähler-Einstein equations built in the previous sections.

Let $(M,g)$ be an ACHE manifold in the sense explained at the end of 
section~\ref{sec:asymptotics}, so that
\[  g - \left(dr^2 + \exp^{2r}\eta^2 + \exp^r\gamma\right) \in  
\cw_{\delta} .\]
From the work done above, one may endow a collar neighbor­hood of 
infin­ity of the form $X^3× ]R,+\infty[$ in $M$ with an 
(integrable at high order) almost complex structure, denoted by $J$. 
This provides us with an approximate Kähler-Einstein metric $\bg$, 
up to the order $O(\exp^{- 3\, r})$ (or even higher, as already noticed)
such that 
\[ \bg - \left(dr^2 + \exp^{2r}\eta^2 + \exp^r\gamma\right) \in  
\cw_{1}. \] 

\begin{lem}\label{lem:gauge}
Let $g$ be an \ache metric and $\bg$ any highly approximate 
Kähler-Einstein metric induced around infinity by the same
CR-structure. Then, for $R$ large enough, there exists
a diffeomorphism $\varphi$ of $[R,+\infty[× X$
inducing the identity at infinity such that $\widetilde{g} =\varphi^*g$
satisfies 
\[ \widetilde{g} - \bg \in \cw_{\delta}  \ \textrm{ and } \ 
\delta^{\bg}\, \varphi^*g + \frac{1}{2}\, d\tr_{\bg} \varphi^*g = 0
\ \textrm{ on } \ [R,+\infty[× X. \]
\end{lem}

\begin{proof}
First, letting $\chi$ be a cut-off function with value $1$ in $B(R')$ 
and value $0$ in $M-B(2R')$, the $C^{k,\alpha}_{\delta}$-norm of 
\[ g- \left( \chi g + (1-\chi) \bg\right)  \] 
(for a given pair $k\in \NM$, $\alpha\in ]0,1[$) 
may be chosen very small if $R'$ is chosen large enough. 
From now on, we denote by $\bg$ the metric $\chi g + (1-\chi) \bg$,
which retains the crucial property of being approximately Kähler-Einstein
around infinity, and (up to enlarging $R'$ again) has strictly negative
Ricci curvature. From \cite[Prop. I.4.6]{Biq00}, there
exists a unique diffeomorphism $\varphi$, approximated at infinity by the
identity up to an element of regularity $C^{k+1}$ and order 
$O(\exp^{-\delta r})$, such that
\begin{equation}\label{eq:gauge}
\delta^{\bg}\, \varphi^*g + \frac{1}{2}\, d\tr_{\bg} \varphi^*g = 0
\end{equation}
and the Lemma is proved.\end{proof}

\medskip

Together with the Einstein equation $\Ric^g + \frac{3}{2}\, g =0$ 
(which is of course preserved by the action of diffeomorphisms), 
$\widetilde{g}=\varphi^*g$ is 
a solution of an elliptic non-linear system of equations, which might be 
written as
\begin{equation}\label{eq:defPhi}
\Phi^{\bg}(\widetilde{g}) \ := \ \Ric^{\widetilde{g}} + 
\,\frac{3}{2}\,\widetilde{g} + 
\left(\delta^{\widetilde{g}}\right)^*
\left(\delta^{\bg}\, \widetilde{g} + \frac{1}{2}\, d\tr_{\bg}\widetilde{g}
\right) \ = \ 0 
\end{equation}

\medskip

The final step of this section is then achieved by showing the

\begin{prop}\label{pro:poids}
The difference $\widetilde{g}-\bg$ lives in the weighted space
\[ \bigcap_{\eps >0 } \cw_{2-\eps} .\]
\end{prop} 

\begin{rem}
As all arguments below will concern behavior at infinity only, one may
consider ({e.g.} by interpolating between the classical complex hyperbolic
metric and $\bg$) that both the reference metric $\bg$ and the \ache metric
$g$ have nonpositive sectional curvature. This has the sole effect that
$g$ (hence $\widetilde{g}$) solves the Einstein equation, or, equivalently,
equation (\ref{eq:defPhi}), up to a compactly supported perturbation only. 
This convention will be in order throughout the remaining parts of this 
paper.
\end{rem}

\medskip

\begin{proof}
The linearization of the map $\Phi^{\bg}$ is computed in 
\cite[formula (I.1.9)]{Biq00} and reads:
\begin{equation}\label{eq:dPhi}
d_{\bg}\Phi^{\bg}(h) = \frac{1}{2}\left(\nabla^{\bg}\right)^*\nabla^{\bg} h
- \Rsym{}^{\bg}(h) 
+ \frac{1}{2}\left( \Ric^{\bg}\circ h + h \circ \Ric^{\bg}
+ 12\, h \right). 
\end{equation}
From the main results of \cite[I.4.B]{Biq00}, the basic 
isomorphism (Proposition I.2.5 of \cite{Biq00}) shows that 
$d_{\bg}\Phi^{\bg}$ is an isomorphism  in weighted Hölder spaces
$C^{k,\alpha}_{\delta}$ for each $(k,\alpha)$, whenever $\bg$ has 
nonpositive curvature
and (in dimension 4) $ 0 < \delta < 2$.  Moreover, since 
$\widetilde{g}-\bg$ may be taken to have small norm in 
$C^{k,\alpha}_1$-topology, one may write
\[ 0 = \Phi^{\bg}(\widetilde{g}) \ = \ \Phi^{\bg}(\bg) + d_{\bg}\Phi^{\bg}
(\widetilde{g}-\bg) + 
P_1(\widetilde{g}-\bg), \]
where $P_1$ is a quadratic term in $g-\bg$. From this we deduce
\begin{equation} d_{\bg}\Phi^{\bg}(\widetilde{g}-\bg) \in C^{k-2,\alpha}_2 
\textrm{ and } \widetilde{g}-\bg \in C^{k,\alpha}_{\delta} . \end{equation}
The isomorphism theorem shows then that $\widetilde{g}-\bg$ lives in 
$C^{k,\alpha}_{2-\eps}$ for each $(k,\alpha)$ and any $\eps >0$. \end{proof}

\bigskip

\section{High-order asymptotic expansion}\label{sec:da}

\bigskip

The goal of this section is to improve the previous asymptotic 
expansion for an \ache metric $g$. In the previous section, we have shown 
that, up to diffeomorphism action, its expansion up to order 
$\exp^{-(2-\varepsilon)r}$ ($\varepsilon>0$) is exactly the same as the
one of the (approximately) Kähler-Einstein metric \(\bg\).
We shall now study what happens at order $\exp^{-2r}$. 

We denote by the same letter $g$ the metric that was denoted by 
$\widetilde{g}$ in the previous section. Hence it satisfies the 
conclusions of Lemma~\ref{lem:gauge} and Proposition~\ref{pro:poids}. 

\medskip

Let $u =g-\bg$. From the previous section
$u$ lives in the weighted space $\cw_{2-\varepsilon}$ for any
(small) $\varepsilon>0$. 
Using the notations
of the previous section, this shows $P_1(u)$ is in 
$\cw_{4-2\varepsilon}$ hence
\[ d_{\bg}\Phi^{\bg}(u) + \Phi^{\bg}(\bg)  \in \cw_{2+\delta},\quad 
\textrm{for some } \ \delta>0. \]
Note also that the metric $\bg$ is given by a series 
in $\exp^{-\frac{k}{2}r}$ ($k\in\NM$) with coefficients smooth in the boundary
variables, so that
$\Phi^{\bg}(\bg)$ may also be taken in $\cw_{2+\delta}$ with
$\delta>0$ (in section \ref{sec:approx}, we took $\delta =\frac{1}{2}$)
and this stands also true for all derivatives of any order
$(\wn_h)^k(\wn_R)^{\ell}\Phi^{\bg}(\bg)$.

It has been pointed out in section \ref{sec:asymptotics} (Lemma 
\ref{lem:commut}) that the metric Laplacian $\nabla^*\nabla$ of the 
$0$-th order metric $g_0 = dr^2 + \exp^{2r}\eta^2 + \exp^r \gamma$ enjoys
nice commutation properties with the (extended) Webster covariant
derivatives $\wn$ in transverse directions. We
shall now use this important piece of information to get a refined
description of the asymptotic behaviour of $u$.
Referring to Lemma \ref{lem:commut} above, we begin by stating:

\begin{lem}\label{lem:com2}
Let $L(h) =\Rsym{}^{\bg}(h) 
- \frac{1}{2}\left( \Ric^{\bg}\circ h 
+ h \circ \Ric^{\bg}+ 12\, h \right)$. Then the operator $d_{\bg}\Phi^{\bg} = 
\frac{1}{2}\Delta_{\bg} - L$ enjoys the 
same commutation properties as the metric Laplacian of $g_0$ with respect to
the Webster covariant derivatives $\widetilde{\nabla}^W_{\partial_r}$,
$\widetilde{\nabla}^W_R$, and $\widetilde{\nabla}^W_h$, $h$ being a 
$\gamma$-unit local section of $H$ on $X$.
\end{lem}

\begin{proof} First of all, $\bg - g_0$, hence $\Delta_{\bg} - \Delta_{g_0}$, 
is in $e^{-r}\cQ_1$. As a result, commutations properties shown
for $\Delta_{g_0}$ remain valid for $\Delta_{\bg}$.
Moreover, the difference between $\frac{1}{2}\Delta_{\bg}$ and 
$d_{\bg}\Phi^{\bg}$ 
is the $0$-th order term $L$ whose coefficients involve the curvature of $\bg$. 
From Corollary \ref{cor-nablaa}, Webster-derivatives of those are 
$O(\exp^{-r})$, whence are elements of $\exp^{-r}\cQ_1 \subset 
\cQ$. \end{proof}

\medskip

We now prove the main result of this section:
 
\medskip

\begin{prop} \label{prop:ke-2r}
There exists a smooth section $k$ of the anti-$J_0$-invariant
symmetric bilinear forms on the contact distribution $H$ on $X$ such
that 
\[ g-\bg - k\exp^{-r} \ \in \ C^{\infty}_{2+\delta} \]
all for some $\delta > 0$. The same is true for all Webster covariant
derivatives of any order.
\end{prop}

Checking decays at infinity, this result amounts to say that $u = g -\bg$ 
is the sum of a smooth term of type $u_{\infty}\exp^{-2r}$ and 
remainder terms of extra decay at infinity (the reader is warned that, 
in the norms associated to $g$ or $\bg$, any term of the type $k\exp^{-r}$
with $k$ defined on $X$ as above 
is indeed of order $\exp^{-2r}$). The main
difficulty of the proof of the Proposition is getting the full 
$C^{\infty}$
regularity on either the leading term $u_{\infty}\exp^{-2r}$ or
the remainder.

\smallskip

\begin{proof} We shall write $P$ for the linear operator 
$d_{\bg}\Phi^{\bg}$ and $f=\Phi^{\bg}(\bg)$, so that 
\[ Pu \ = \ f + P_1 u \ \in \ C^{\infty}_{2+\delta} \ \textrm{ for } 
\delta>0 ,\]
where $P_1$ is the non-linear part of the Einstein equations
(with its extra gauge fixing terms). Moreover 
$(\wn_R)^k (\wn_h)^{\ell}f$ belongs
to $C^{\infty}_{2+\delta}$ for any $k,\ell$, too.
 
\begin{sulem}
Let $w\in C^{\infty}_{\alpha}$ for some $\alpha >0$. 
If $Pw \in C^{\infty}_{2+\delta}$ then, for any $Q\in\cQ$, 
\[  (\p_r+2)w \in C^{0}_{2+\delta}
\ \textrm{ and } \  Qw \in C^{\infty}_{2+\delta}. \]
\end{sulem}

{\flushleft\it Proof of the Sub-lemma}. -- First of all, elementary
weight considerations as above show that $w$ lies in 
$C^{\infty}_{2-\varepsilon}$ for every $\varepsilon>0$.
From the commutations properties proved in Lemma \ref{lem:commut}, one has
\[ P(\wn_R w) \ = \ \wn_R (Pw) + [P,\wn_R] w \  \in C^{\infty}_{1+\delta}
\]
for some $\delta >0$. Indeed, $\wn_R (Pw)$ obviously lies in 
$C^{\infty}_{1+\delta}$ whereas the bracket term preserves weights as it 
belongs to $\cQ$, a subset of the weight-preserving operators.
We then get that
\[ \wn_R w \in C^{\infty}_{1+\delta} .\]
Moreover, it exists $Q\in\cQ$ such that, for every $\gamma$-unit
$h$ in $H$ 
\[ P(\wn_{h} w) \ = \ \wn_{h} (Pw) - 2 \wn_{\exp^{-\frac{r}{2}}Jh}
\wn_{\exp^{-\frac{r}{2}}R} w  + Qw \ ,\]
so that $P(\wn_{h} w)$, hence $\wn_{h}w$, belongs to $C^{\infty}_{3/2+\delta}$
(for both arguments, we recall that the critical weights of $P$ are 
$0$ and $2$). Obviously, this implies that $Qw$ lies in $C^{\infty}_{2+\delta}$
for each $Q$ in $\cQ$.

Using this control on the transverse derivatives, the fact that 
$Pw\in \cw_{2+\delta}$ now translates into
\[
q(w) \ = 
\ \left( -\partial_r^2 - 2\partial_r + A \right) w 
\ \in \ \cw_{2+\delta}
\]
where $A$ is a linear operator on symmetric bilinear forms, obtained
as the dominant term in the asymptotic expansion of the Laplace operator
and the $0$-th order terms in (\ref{eq:dPhi}). This follows
from the work done in \cite[Section I.2]{Biq00}. The exact expression
of $A$ will be of no concern to our purposes; the following information 
is nonetheless crucial: from \cite[Section I.4.B]{Biq00}, we know that 
the smallest eigenvalue of $A$ is equal to $0$, and the associated 
eigen-subbundle is the 
bundle of anti-$J_0$-invariant symmetric bilinear forms on the contact 
distribution $H$, {\sl i.e.} the quadratic forms $k$ on $TX$ such that 
\[ k(R,\cdot)=0, \ \ \  
k(J_0\cdot,J_0\cdot) \, = \, - \, k(\cdot,\cdot). \]
The highest critical weights of $P$ are then $2$ on this bundle and
are larger on its orthogonal complement.
This has two consequences: first of all, we get that all components
of $w$ that are (pointwise) orthogonal to the eigen-subbundle are
elements of $C^{\infty}_{2+\delta}$ (for some $\delta>0$), whereas
elementary ordinary differential equations analysis applied to the
remaining component shows that there is a symmetric field $w_\infty$ on 
$X$ as above, such that
\begin{equation}\label{eq:asymp-precise}
w - w_\infty \exp^{-2r} \in C^0_{2+\delta}\ .
\end{equation}
More precisely, let us project the solution $w\in \cw_{2-\varepsilon}$ 
and the control $q(w)\in\cw_{2+\delta}$, on the eigenbundle of 
anti-$J_0$-invariant symmetric bilinear forms on $H$; using for sake of 
simplicity the same letters for the projections, we 
then have $q(w)=-\p_r^2 w-2\p_r w$, so that we may write 
$w=w_\infty \exp^{-2r}+w'$, with
\begin{equation}\label{eq:wprime}
 w' = e^{-2r} \int_{\infty}^r e^{2t} \int_t^{\infty} q(w)\, ds
\end{equation}
obtained by integration along each ray $\{x\}\times [R_0,+\infty)$.
Then $w'$ and $\p_r^k w'$ (for any $k$) are easily seen to be sections
of $C^0_{2+\delta}$. The control on $(\p_r+2)w=(\p_r+2)w'$ now
follows. \qed

\medskip

{\flushleft\it Proof of Proposition \ref{prop:ke-2r} \upn{(}continued\upn{)}}.
-- It results from the previous analysis (applied to $w=u= g-\bg$)
that $u=\exp^{-2r}u_\infty + u'$, but the tangential regularity on $u_\infty$
and $u'$ is not yet known. We now prove by induction on 
$(\ell+\ell',\ell)$
the following assertion: let $h_1,\dots,h_{\ell}$ be sections of $H$ 
on $X$, then, for any $Q\in\cQ$
$$Q\wn_{h_1}\cdots\wn_{h_{\ell}}(\wn_R)^{\ell'} u\in \cw_{2+\delta}.$$
Observe that the order of the derivations here is not important, since 
they
commute up to introducing terms with less derivatives.

Before explaining the induction, let us observe how Proposition 
\ref{prop:ke-2r} is a
consequence of the assertion. Indeed, it implies that the tangential
derivatives of $Pu$ are controlled, so that one has, for each $\ell$,
$\wn_{h_1}\dots\wn_{h_{\ell}}q(w)\in\cw_{2+\delta}$,
so by \eqref{eq:wprime} one 
gets $\wn_{h_1}\cdots\wn_{h_{\ell}}(\wn_R)^{\ell'}u'\in\cw_{2+\delta}$.
The regularity of $u_\infty$ follows.

We now come back to the induction. The case $\ell+\ell'=0$ is 
the Sub-lemma above.
So now fix $(\ell,\ell')$ and suppose that the result is true for all 
$(\ell_1,\ell_1')$
such that $\ell_1+\ell_1'<\ell+\ell'$, or $\ell_1+\ell_1'=\ell+\ell'$ 
and $\ell_1<\ell$.
Therefore, we control $\ell+\ell'$ transverse derivatives for $u$ in
$\cw_{1+\delta}$.
As $P_1$ is quadratic-or-more with
coefficients given by smooth functions on the boundary at infinity, we
certainly have also a control of it
in $\cw_{2+\delta}$ for $\ell+\ell'$ transverse
derivatives of $P_1u$. As $Pu = f + P_1u$ and $f$ is controlled in the
best way one can hope, it follows therefore that all $\ell+\ell'$ 
transverse derivatives $Pu$ also live in $\cw_{2+\delta}$. 
Letting $\xi_i=h_i$ for $i\leq \ell$, $\xi_i=R$ for $i>\ell$, 
and $(\wn_{\xi})^{\ell+\ell'}$ for 
$\wn_{\xi_{1}}\cdots\wn_{\xi_{\ell+\ell'}}$, we can write
\begin{equation}\label{eq:recurrence}
P\left((\wn_{\xi})^{\ell+\ell'}u\right)
=(\wn_{\xi})^{\ell+\ell'}Pu + 
 \sum_{i=1}^{\ell+\ell'} \wn_{\xi_1}\cdots\wn_{\xi_{i-1}}[P,\wn_{\xi_{i}}]
                  \wn_{\xi_{i+1}}\cdots\wn_{\xi_{\ell+\ell'}}u .
\end{equation}
We would like to argue that all the terms in the sum are in 
$\cw_{2+\delta}$. For this we have to distinguish two cases from 
the commutation rules for $P$ in Lemma \ref{lem:com2}, for  
bracket terms in (\ref{eq:recurrence}), 

\smallskip

(i) if $i>\ell$, then
in this case $\xi_i=R$ and $[P,\wn_{R}]\in\cQ$, so that the corresponding 
bracket term lies in $\cw_{2+\delta}$ by the induction 
hypothesis for $\ell+\ell'-1$;

\smallskip

(ii) if $i\leq \ell$, then in this case $\xi_i=h_i$ and we know that
\[ [P,\wn_{h_{i}}]=-2\,\nabla_{\exp^{-\frac{r}{2}}Jh_i}
\wn_{\exp^{-\frac{r}{2}}R} \ \textrm{ mod } \ \cQ  \]
so that, modulo $\cQ$, the bracket term becomes
\[ 
-2\,\nabla_{\exp^{-\frac{r}{2}}Jh_i}
\left(\wn_{\xi_1} \cdots \wn_{\xi_{i-1}}
\wn_{\exp^{-\frac{r}{2}}R}\wn_{\xi_{i+1}}\cdots\wn_{\xi_{\ell+\ell'}} u 
\right) \]
and the induction hypothesis with $\ell+\ell'$ unchanged but a smaller 
$\ell$ shows that this term has fast decay (remember that
$\nabla_{\exp^{-\frac{r}{2}}Jh_i}$ is in $\cQ$)
and therefore the term under consideration in the right-hand side of 
(\ref{eq:recurrence}) above is again in $\cw_{2+\delta}$.

\smallskip

We get at the end that each bracket term in (\ref{eq:recurrence}) is in 
$\cw_{2+\delta}$, and it follows that
$P\wn_{h_{1}}\cdots\wn_{h_{\ell}}u$ also is in $\cw_{2+\delta}$.
Applying the Sub-lemma to $(\wn_{\xi})^{\ell+\ell'}u$ yields the desired
induction statement for $(\ell+\ell',\ell)$.
\end{proof}

\medskip

When injecting the precise asymptotic expansion of the Kähler-Einstein
metric $\bg$, the final output is

\begin{cor}[asymptotic expansion]\label{asymptotic-expansion}
Let $(M,g)$ be any asymptotically complex hyperbolic Einstein manifold 
of dimension $4$ with CR-structure at infinity $(\eta,H,J_0)$ on a
$3$-manifold $X$, and $\bg$ the approximate Kähler-Einstein metric
determined by the CR-structure. Then, 
up to a $C^{k+1}$-diffeomorphism inducing the identity at infinity, $g$ 
has an asymptotic expansion of the form
\begin{equation*}\label{eq:DAache}
g = \bg +  k\exp^{-r} + \text{\rm\ lower order terms}
\end{equation*}
where $k$ is a anti-$J_0$-invariant quadratic form on the contact 
distribution $H$. The value of $k$ is formally undetermined. 
\end{cor}

Together with the result of Corollary \ref{cor:omega-explicit}, this 
gives 
the asymptotic expansion promised in remark \ref{remann}.
Once again, the reader must take care of the fact that the term
$k\exp^{-r}$ grows like $\exp^{-2r}$.

\bigskip

\section{Convergence of the integral}\label{sec:weyl}

\bigskip

In this section, we prove the convergence of our integral for any ACHE
metric. We will prove the convergence of two terms, the first
involving only the antiselfdual Weyl tensor, the second one involving
the selfdual Weyl tensor and the scalar curvature. In the
Kähler-Einstein case, only the first term occurs and the proof
simplifies considerably, as shown by the following:

\begin{lem}
If $(M,g)$ is an ACH Einstein manifold, then $|W^-|=O(\exp^{-\delta r})$
for any $\delta<2$. In particular,
$$ \int |W^-|^2 < \infty. $$
\end{lem}

\smallskip

\begin{rem}
The behavior of $W^-$ will be precised further in the course of the proof,
see formula (\ref{behWm}).
We will show that the highest-order of the asymptotic 
expansion of $W^-$ occurs at a decay $\exp^{-2r}$ and that this term 
belongs to a specific sub-bundle of the bundle of Weyl curvature tensors
(notice that the existence of a an asymptotic expansion up to order 
$\exp^{-\frac{5}{2}r}$ at least flows from the existence of
the asymptotic expansion of the metric given in Corollary 
\ref{asymptotic-expansion}). As in the previous section, there are no 
``logarithmic terms'' (or, rather $r^k\exp{-2r}$-terms in our context), 
as these would be formally determined. As the general CR case cannot
be distinguished locally from the embeddable CR case, these terms would 
necessarily appear in the development of the K\"ahler-Einstein solution $\bg$.
But the critical weight of the complex Monge-Amp\`ere is higher (it is
equal to $3$) \cite{Fef76} and this prevents the appearance of any such
terms at order $2$. 

Note moreover that a semi-explicit form of the
highest-order term in $W^-$ in the Kähler-Einstein case will be given
in Proposition \ref{weyl-cartan} after the current lemma.
\end{rem}

\smallskip

\begin{proof}
We use the fact that, for an Einstein metric, the Weyl tensor is
harmonic as a 2-form with values in the endomorphism of $TM$. Therefore
$$
\left( (d^\nabla)^*d^\nabla+d^\nabla(d^\nabla)^* \right)W^-=0.
$$
By \cite[proposition I.3.5]{Biq00}, any such harmonic form which is
$O(\exp^{-(\delta_-+\eps)r})$ must be
$O(\exp^{-\delta r})$ for any $\delta<\delta_+$, where
$$\delta_{±}=1±\sqrt{1+\lambda}$$
are the critical weights of the associated indicial operator and 
$\lambda$ is the smallest eigenvalue of the ``zero order terms'' of
the operator $(d^\nabla)^*d^\nabla+d^\nabla(d^\nabla)^*$ at infinity
(see \cite{Biq00} or \cite{DeHe01} for details).

This $\lambda$ does not depend on the particular conformal infinity
$\gamma$, and can be calculated for the model $\chii$.
We now restrict to the case of a harmonic antiselfdual form $w$ with
values in $\Omega^2_-$ on $\chii$. There is a Weitzenböck formula
expressing the difference between the Laplacians $\nabla^*\nabla$ and
$(d^\nabla)^*d^\nabla+d^\nabla(d^\nabla)^*$ as an algebraic operator
involving the curvatures of the manifold and of the bundle, see
\cite[theorem 3.10]{BouLaw81}; in our case, because $W^-(\chii)=0$,
the term involving the curvature of the basis reduces to
$\frac{\Scal}{3}$, so we get
$$
\left( (d^\nabla)^*d^\nabla+d^\nabla(d^\nabla)^* \right) w =
\nabla^* \nabla w + \frac{\Scal}{3}w + \mathcal{R}(w),
$$
where $\mathcal{R}(w)$ is the 2-form with values in the endomorphisms of
$TM$ given by
\begin{equation}\label{def-RR}
\mathcal{R}(w)_{X,Y}=\sum_1^4 \left(
  [R^\nabla_{e_j,X},w_{e_j,Y}]-[R^\nabla_{e_j,Y},w_{e_j,X}]
                              \right).
\end{equation}
When $w$ takes its values in $\Omega^2_-$, only the antiselfdual part
of the curvatures $R^\nabla_{e_j,X}$ may give a nonzero result in the
equality, which means that only the scalar curvature and $W^-$ are
involved. Since again $W^-(\chii)=0$, this means that we can calculate
$\mathcal{R}$ using the constant sectional curvature tensor 
$$ R_{X,Y}=-\frac{\Scal}{12}X\land Y. $$
Now an easy explicit calculation using \eqref{def-RR} gives actually,
for a section $w$ of $\Sym_0^2\Omega^2_-$,
$$ \mathcal{R}(w)=\frac{\Scal}{6}w. $$

We now look at the Laplacian
$\nabla^*\nabla$. To state the result of the computation for this term, 
we first need to decompose the bundle $\Omega^2_-$ over the
boundary: we have the local orthonormal basis
$$
(e_1=\partial_r,e_2=Je_1=\exp^{-r}R,e_3=\exp^{-\frac{r}{2}}h,
e_4=Je_3=\exp^{-\frac{r}{2}}Jh),
$$
and this induces a local decomposition
$$
\Omega^2_- = \RM \oplus \CM
$$
with $\RM$ generated by $e^1e^2-e^3e^4$ (all wedge products are suppressed
in the rest of this section), and $\CM$ generated by
$e^1e^3+e^2e^4$ and $e^1e^4-e^2e^3$; more intrinsically, $\CM$ is
isomorphic to the contact distribution $H$ of the boundary (or more
precisely its dual $H^*$).
Then inside $\Sym_0^2\Omega^2_-$ we find the real 2-dimensional
subbundle $\Sym_0^2\CM=\Sym_0^2 H^*$, which consists of
anti-$\CM$-linear endomorphisms of $H^*$, or equivalently of $H$:
these are the infinitesimal deformations of the complex structure on
$H$.

\smallskip 

\begin{notn} We will denote by $\cJ$ the bundle of infinitesimal 
deformations of the complex structure in $H$, seen as a subbundle of 
$\Sym_0^2\Omega^2_-$.
\end{notn}

\smallskip

\begin{claim}
The smallest eigenvalue of the zero order term of $\nabla^*\nabla$
acting on $\Sym_0^2\Omega^2_-$ is $-\frac{\Scal}{2}$.
The corresponding eigenspace is $\cJ$.
\end{claim}

\smallskip

First admit the claim: then the smallest eigenvalue of the zero order
term of the Laplacian $(d^\nabla)^*d^\nabla+d^\nabla(d^\nabla)^*$ is
finally
$$ \lambda=\frac{\Scal}{3}+\frac{\Scal}{6}-\frac{\Scal}{2}=0. $$
Therefore $\delta_+=2$ and $\delta_-=0$, which means that a harmonic
section of $\Sym_0^2\Omega^2_-$ which goes to 0 at infinity must
actually decay as $O(\exp^{-2r})$, and the term of order~2 lies in the
subbundle $\cJ$. Apply this to the tensor $W^-$ of an ACHE metric: for
any ACH metric, we certainly have $|W^-|=O(\exp^{-\eps r})$, and therefore
\begin{equation}\label{behWm}
W^- = W^-_2 \exp^{-2r} + O(\exp^{-(2+\eps)r}),\quad W^-_2\in\cJ.
\end{equation}

It remains to prove the claim. One can represent the hyperbolic space
as a quotient $U_{1,2}/U_1U_2$.
By \cite[I,(2.12)]{Biq00}, one has to calculate the action of a
partial Casimir operator on the representation $\Sym^2_0\Omega^2_-$ of
the group $U_1U_2$. Analogous calculations are performed in
\cite[I.4.B]{Biq00} or in \cite{DeHe01}, and the details of the present 
case are left to the reader.
\end{proof}

\medskip

The proof of the lemma is valid also for the approximate Kähler-Einstein 
metric $\bar{g}$.
Therefore, we also have
$$
W^{-}(\bar{g})=W^{-}_{2}(\bar{g}) \exp^{-2r} + O(\exp^{-(2+\varepsilon)r}), \quad W^-_2(\bar{g})\in\cJ.
$$

Remark that the filling complex structure in $X× ]R,+\infty)$ and the formally
determined part of the Kähler-Einstein metric $\bar{g}$
depend only on the CR structure underlying the pseudohermitian one.
Therefore, $W^{-}_{2}(\bar{g})$ is a CR local invariant.
From the formula for the metric $\bar{g}$ in Corollary
\ref{cor:omega-explicit}, we see that it depends only on four
(horizontal) derivatives of the CR structure, and linearly only of the
fourth derivatives.

By the classical work on CR geometry of Chern and Moser \cite{ChMo74}
(actually going back in that case to Cartan), the only such invariant 
is the Cartan curvature tensor $Q$, which can be seen also as a
section of the bundle $\cJ$. Therefore we get:

\begin{prop}\label{weyl-cartan} Let $Q$ be the Cartan tensor of the 
CR structure, seen as an anti-selfdual Weyl-type tensor. Then
$$W^-_{2}(\bar{g}) = a\, Q\,+ O(\exp^{-(2+\varepsilon)r}).$$
\end{prop}
The constant $a$ can be determined from section \ref{sec:nu-variation}.
Of course one should also be able to calculate directly the curvature
tensor of $\bar{g}$ from the formula of corollary \ref{cor:omega-explicit}:
this would give another proof of the proposition, and also the value
of the constant $a$.

\medskip
We now turn to the $W^{+}$-term in the integral.
\begin{lem}
If $(M,g)$ is an ACH Einstein metric, then
$$ \int |W^+|^2 - \frac{1}{24} \Scal^2 < \infty.$$
\end{lem}
\begin{proof}
In the Kähler-Einstein case, one has identically
$$ |W^+|^2 = \frac{1}{24} \Scal^2 $$
and there is nothing to prove. In the general case, we will use the
fact, from corollary \ref{asymptotic-expansion}, that an ACH Einstein
metric differs from an asymptotically Kähler-Einstein metric $\bg$ up
to a term of order $\exp^{-2r}$, namely 
$$ g=\bg + k \exp^{-r} + O(\exp^{-(2+\eps)r}). $$
Because the volume grows like $\exp^{2r}$, we need to look only at
terms of order at most $\exp^{-2r}$. In particular, we can neglect the
fact that $\bg$ is Kähler-Einstein only up to order $\exp^{-3r}$ for 
example.

Because the term $k \exp^{-r}$ is only a small perturbation near
infinity, the difference between the Weyl tensors of $g$ and $\bg$
depends on the difference $g-\bg$ by a linear part (the differential
$d_{\bg}W$) and a quadratic part which has a stronger decay. Therefore,
$$
W^+(g)=W^+(\bg)+d_{\bg}W^+(k \exp^{-r})+O(\exp^{-(2+\eps)r}),
$$
and the lemma will be proved if we are able to prove that $d_{\bg}|W^+|^2(k
\exp^{-r})$ is $O(\exp^{-(2+\eps)r})$. As the metric $\bg$
is ACH, it is enough to compute
$d_{\bg}W^{+}$ when $\bg$ is the standard complex hyperbolic metric where
of $\chii$, the difference with the
general terms giving only lower order terms (this is a consequence
of Remark \ref{rem:reexpressed}).

The differential of the Weyl tensor is well-known and we extract the 
following useful facts from \cite[2.5]{Gau93}: for a path of metrics 
$g_{t}$, the symmetric endomorphism 
$u_{t}=(g_{0}^{-1}g_{t})^{\frac{1}{2}}$ sends 
the metric $g_{0}$ to the metric $g_{t}$. If we consider the curvatures 
$R=R(g_{t})$ as sections of $\Omega ^{2}\otimes \Omega ^{2}$, then
\begin{equation}\label{eqdR1}
\dR = \mathcal{R} + \frac{1}{2}\ad_{\dg}R,
\end{equation}
where $\mathcal{R}$ is the derivative of $u_{t}^{-1}R_{t}$, and
$$
(\ad_{\dg}R)(X_{1},\cdots ,X_{4})=
-\sum _{1}^{4}R(X_{1},\cdots ,\dg(X_{i}),\cdots ,X_{4}).
$$
As we are interested only in $d_{\bg}|W^+|^2(k\exp^{-r})$, we will only 
need the value of the projection 
$w^+=\pi_{\cW^+}\mathcal{R}$ of $\mathcal{R}$ on
$\cW^+=\Sym_0^2\Omega^2_+$ (as
$d_{\bg}|W^+|^2 = d_{\bg}|u_t^{-1}W^+(g_t)|^2$, the other terms do 
not contribute).

The variations of the Weyl tensors depend only on the deformation of 
the conformal structure. This is parameterized by the deformation of 
$\Omega^2_-$, given by some $u\in\Omega^2_-\otimes\Omega^2_+$, and
from \cite[(2.5.8)-(2.5.12)]{Gau93},
$$
w^+ = -\frac{1}{2} \Pi d^\nabla_+ (d^\nabla)^* u,
$$
where $(d^\nabla)^* u$ belongs to $\Omega^1\otimes\Omega^2_+$, the operator
$d^\nabla_+$ is the projection of $d^\nabla$ on selfdual 2-forms, and
$\Pi$ is the projection
$\Omega^2_+\otimes\Omega^2_+\to\Sym^2_0\Omega^2_+$.

We can now get explicitely this differential. Coming back to the
calculation of the Levi-Civita connection in lemma \ref{lem-nabla}, we
use again the orthonormal frame
\begin{equation}\label{frame}
(e_1,e_2,e_3,e_4)=(\partial_r,\exp^{-r}R,\exp^{-r/2}h,\exp^{-r/2}Jh).
\end{equation}
We will use also the following basis of selfdual and antiselfdual 2-forms:
\begin{equation}\label{omegai}\begin{split}
\omega^1_±&=e^1e^2± e^3e^4, \\
\omega^2_±&=e^1e^3\mp e^2e^4, \\
\omega^3_±&=e^1e^4± e^2e^3.\end{split}
\end{equation}
The tensor $k$ is given in the basis $(e_3,e_4)$ by the tracefree
symmetric matrix 
\begin{equation}\label{defk}
k=2\exp^{-r}\begin{pmatrix}a&b\\b&-a\end{pmatrix}.
\end{equation}
Infinitesimally, the orthonormal basis $(e^j)$ is transformed into
$$
(e^1,e^2,e^3+\exp^{-2r}(ae^3+be^4),e^4+\exp^{-2r}(be^3-ae^4)).
$$
Let us calculate the modification of $\Omega^2_-$: the form
$\omega^2_-$ becomes
$$
e^1(e^3+\exp^{-2r}(ae^3+be^4))+e^2(e^4+\exp^{-2r}(be^3-ae^4))=
\omega^2_-+\exp^{-2r}(a\omega^2_++b\omega^3_+),
$$
and similarly $\omega^3_-$ is transformed into
$ \omega^3_--\exp^{-2r}(a\omega^3_++b\omega^2_+) $
and $\omega^1_-$ is unchanged; this means that the deformation
$k\exp^{-r}$ of the metric translates into the tensor
$$
u=a\exp^{-2r}(\omega^2_-\omega^2_+-\omega^3_-\omega^3_+)
 +b\exp^{-2r}(\omega^2_-\omega^3_++\omega^3_-\omega^2_+).
$$
From lemma \ref{lem-nabla}, an easy calculation gives us the higher
order terms of the covariant derivatives of the $\omega_2^{±}$:
\begin{align*}
\nabla \omega^2_+ &= \frac{3}{2}e^2\omega^3_+, &
\nabla \omega^2_- &=-\frac{1}{2}e^2\omega^3_-+e^4\omega^1_-, \\
\nabla \omega^3_+ &=-\frac{3}{2}e^2\omega^2_+, &
\nabla \omega^3_- &= \frac{1}{2}e^2\omega^2_--e^3\omega^1_-.
\end{align*}
Using these formulas, we can calculate $(d^\nabla)^*u=-\tr \nabla u$,
keeping in mind that we can neglect derivatives along the boundary,
because these give lower order terms: we get successively
$$
\nabla(\omega^2_-\omega^2_+-\omega^3_-\omega^3_+)=
  e^2(-\omega^2_-\omega^3_+-\omega^3_-\omega^2_+)
 +e^4\omega^1_-\omega^2_++e^3\omega^1_-\omega^3_+,
$$
which leads to
$$
(d^\nabla)^*(\omega^2_-\omega^2_+-\omega^3_-\omega^3_+)=0,
$$
and similarly,
$$
(d^\nabla)^*(\omega^2_-\omega^3_++\omega^3_-\omega^2_+)=0.
$$
Therefore it remains, because of the differentiation of $\exp^{-2r}$
with respect to $r$,
\begin{align*}
(d^\nabla)^*u&=2\tr e^1u \\
&= 2\exp^{-2r}\left(
     a(e^3\omega^2_+-e^4\omega^3_+)+b(e^3\omega^3_++e^4\omega^2_+)
              \right).
\end{align*}
In a similar way, it is straightforward to calculate
$d^\nabla_+(d^\nabla)^*u$. Still restricting to higher order terms,
we give only the result:
$$
d^\nabla_+(d^\nabla)^*u = 3\exp^{-2r}\left(
  a \big(-(\omega^2_+)^2+(\omega^3_+)^2\big)
 -b \big(\omega^2_+\omega^3_++\omega^3_+\omega^2_+\big)
                                     \right).
$$
The important fact here is that $d^\nabla_+(d^\nabla)^*u$ lies in the
orthogonal of
\begin{equation}\label{eq-W+}
W^+ = \frac{\Scal}{6} (\omega^1_+)^2
    - \frac{\Scal}{12} \big((\omega^2_+)^2+(\omega^3_+)^2\big)\ ,
\end{equation}
and from this fact we finally deduce that
$$
\bg\big(W^+,d_{\bg}W^+(k\exp^{-r})\big)=O(\exp^{(-2+\eps)r})
$$
and this ends the proof. \end{proof}

\bigskip

\section{The renormalized integral and the invariant at infinity}
\label{sec:comput}

\bigskip

In the previous sections,  we showed that the metric $g$ differs
from the approximate Kähler-Einstein metric $\bg$ by a term 
of type $k\exp^{-r}$ plus lower order terms living in 
$\cw_{2+\varepsilon}$ ($\varepsilon>0$).
Moreover, $k$ is a section of the bundle of quadratic forms on the 
contact distribution $H$ that are $J_0$-anti-invariant, {\sl i.e.} 
\begin{equation}\label{eq:kJinv}
k(J_0\cdot,J_0 \cdot) \ = \ - \ k(\cdot,\cdot).
\end{equation}
This was the key step to prove that the integral
\begin{equation}\label{eq:intfinale}
\frac{1}{8\pi^2}\,\int_M \left( 3\, |W^-|^2 - |W^+|^2 + 
\frac{1}{24}\,\Scal^2 \right)
\end{equation}
converges for any \ache metric $g$. 

The goal of the present section is to make a few steps towards its 
computation, in terms of a new invariant of the CR structure at
infinity, which will be defined below for any CR manifold (even without 
any filling by an Einstein metric).

\medskip

\subsection*{Boundary terms of characteristic classes}  It will be 
convenient to rewrite the integral (\ref{eq:intfinale}) in terms 
of topological invariants together with boundary contributions. 
In order to do so, we shall
need the formula relating the Euler characteristic and signature of
a compact domain with boundary with the expected interior integral
and local and non-local contributions of the boundary. 

\medskip

In the following
formulas, $R$ will always denote the curvature of the $4$-dimen­sional
manifold and, if $D$ is a bounded domain in $M$, the second fundamental 
form of $\p D$ will always be defined as $\II = \nabla\mathrm{n}$, $\mathrm{n}$ being 
the outer unit normal, and  seen here either as a vector-valued $1$-form 
(rather than as an endomorphism) or as a
quadratic form. For $1$-forms $\alpha$ and $\beta$, we let 
$\alpha\land\beta = \alpha\otimes\beta -\beta\otimes\alpha$; if $\beta$
is a $2$-form, 
$\alpha\land\beta (X,Y,Z) = \alpha (X)\beta(Y,Z) +\alpha (Y)\beta(Z,X) 
+\alpha (Z)\beta(X,Y)$. We define for a tensor $F$ in $\otimes^3T^*M$, 
\begin{equation}\label{eq:csum}
\csum(F)(X,Y,Z) = F(X,Y,Z) + F(Y,Z,X) + F(Z,X,Y), \end{equation} 
hence $\alpha\land\beta = \csum(\alpha\otimes\beta)$. If $\mu$ and 
$\nu$ are forms with values in bundles $E$ and $F$, we decide
that $\mu\land\nu$ is the obviously defined form with values in 
$E\otimes F$. Last, if $\rho\otimes\sigma$ is a $3$-form with values in 
$\otimes^3 TM$, we define
\begin{equation}\label{eq:defT}
\tv(\rho\otimes\sigma) = \la d\!\vol_{\p D}, \sigma \ra \,\rho .
\end{equation}
The desired formula for the Euler characteristic reads: 
\begin{equation}\label{eq:chi-bord}\begin{split}
\chi(D) \ = \ & \frac{1}{8\pi^2} \int_D \left( |W|^2 - 
\frac{1}{2}|\Ric_0|^2 + \frac{1}{24}\Scal^2\right) d\!\vol_D \\ 
& \ + \frac{1}{12\pi^2} \int_{\partial D} \tv( \II\land 
\II\land \II ) - \frac{1}{4\pi^2} \int_{\partial D} \tv ( \II \land R )
\end{split}\end{equation}
where, in the right-hand side, curvature or second forms must be seen as 
$1$- or $2$-forms with values in vectors or $2$-vectors. The signature 
formula includes a non-local boundary contribution $\eta(\p D)$,
known as the $\eta$-invariant:
\begin{equation}\label{eq:tau-bord}\begin{split}
\tau(D) = 
\frac{1}{12\pi^2}\int_D (|W^+|^2 - |W^-|^2)d\!\vol_D + \frac{1}{12\pi^2}
\int_{\partial D} \csum( \II(., R(.,.)\mathrm{n})) + \eta (\partial D) 
\end{split}\end{equation}
where curvature is seen here as a $2$-form with values in endomorphisms and
the second form as a quadratic form.

\medskip

\subsection*{Definition of the invariant at infinity $\nu(X)$} 
From the formulas recalled above, 
our integral in a compact domain $D_r$, whose boundary is the distance 
sphere
$S_r = \{r\}× X$ in the bulk manifold $M$, can be written as the 
sum of a topological contribution $\chi(D_r) - 3 \tau(D_r)$ plus
boundary integrals. As $r$ goes to infinity, the characteristic
numbers remain constant (when $r$ is large enough) and our task 
will now be to study the sum
\begin{equation}\label{eq:intbord}\begin{split}
\nu^g(r)=- \frac{1}{12\pi^2} \int_{S_r} & \tv( \II\land \II\land \II ) 
+ \frac{1}{4\pi^2} \int_{S_r} \tv ( \II \land R )\\
& + \frac{1}{4\pi^2} \int_{S_r} \csum( \II(., R(.,.)\mathrm{n})) 
+ 3\, \eta (S_r,g) 
\end{split}\end{equation}
as the sphere $S_r$ grows to infinity. Recall the tensor $R$ in the above
formula is the curvature of the bulk manifold and $\mathrm{n}$ is the
unit normal vector of slices of constant $r$.
We know from the previous section that this term converges, but we
will show that the limit depends only on the CR-structure on the
boundary. More precisely:

\smallskip

\begin{theo}\label{theo:nu}
If $\bg$ is the asymptotically Kähler-Einstein metric associated to
CR-structure on $\partial_\infty M=X$, then one has
$$ \lim_{r\to\infty} \nu^g(r)-\nu^{\bg}(r)=0. $$
In particular,
$$ \nu(X)=\lim_{r\to\infty}\nu^g(r)
      =\lim_{r\to\infty}\nu^{\bg}(r) $$
and the limit is an invariant which depends only on the CR-structure on 
the boundary at infinity.
\end{theo}

This leads us to the general definition of the invariant $\nu$. Indeed,
let $X$ be an arbitrary CR 
manifold of dimension $3$ and let $\overline{M} = X× ]R,+\infty)$ for
some $R>0$. The work done in section \ref{sec:weyl} is independent of
topology, or even completeness, of the filling Einstein metric, hence
it shows that, 
if $\bg$ is the approximate Kähler-Einstein metric (up to order 
$\frac{5}{2}$) on $\overline{M}$
constructed in the first sections, the quantity $\nu^{\bar{g}}(r)$
converges as the sphere $S_r$ grows to infinity, and the limit is a 
CR invariant of $X$. Hence we can define:

\begin{de}\label{def:nugeneral}
Let $X$ be an arbitrary compact, pseudoconvex CR manifold of
dimension $3$. Then 
$$ \nu(X) \ = \lim_{r\to\infty}\nu^{\bg}(r) \ ,  $$
defined with the approximate K\"ahler-Einstein metric $\bg$ of Theorem
\ref{theo:omegaKE}, is an invariant of the CR structure.
\end{de}

\medskip

Now we can state, as a consequence, the following result, proving
Theorem \ref{btheo}.

\begin{cor}\label{theo:dependance}
For any \ache metric $g$, inducing a CR-structure on the boundary at
infinity $\partial_\infty M=X$, one has 
\[
\frac{1}{8\pi^2}\,\int_M \left( 3\, |W^-|^2 - |W^+|^2 + 
\frac{1}{24}\,\Scal^2 \right) = \chi(M) - 3 \tau(M) +
\nu(X) .
\]
\end{cor}

\medskip

\subsection*{Proof of theorem \ref{theo:nu}}
If we denote by $B(g)$ the local boundary integrand in formula 
(\ref{eq:intbord}) computed with the metric $g$, we can write, 
for $r$ large enough,
\begin{equation}\label{eq:bbord}\begin{split}
\int_{S_r} B(g) + 3\,\eta(S_r,g) \ = \ & \int_{S_r} B(\bg) \ 
+ \ 3\,\eta(S_r,\bg)\\
& \ + \ \int_{S_r} (d_{\bg}B)(g-\bg) \ + \ 3 (d_{\bg}\eta)(g-\bg) \\ 
& \ + \ \int_{S_r} {\mathcal Q}(g-\bg), \end{split}
\end{equation}
where, as usual, ${\mathcal Q}$ is a quadratic term, {\sl i.e.} there
exists a constant $C>0$ and some $\epsilon>0$ such that, for any $r$ 
large enough
\begin{equation}\label{eq:controleQ}
{\mathcal Q} (g-\bg) \ \leq \ C\,\exp^{-(2+\epsilon)r}\ .
\end{equation}
From the work by Burns and Epstein already quoted, the first two terms
in the right hand side of (\ref{eq:bbord})
converge as they depend only on the (approximately)
Kähler-Einstein metric $\bg$. Moreover, as $|d\!\vol_{\bg}| = 
O(\exp^{-2r})$, ${\mathcal Q}(g-\bg)$ has integrals converging to
zero, and we are now reduced to show that, under the above assumptions 
for $g$ and $\bg$,
\begin{equation}\label{eq:limite}
\lim_{r\to\infty} \ \int_{S_r} (d_{\bg}B)(g-\bg) \ 
+ \ 3 (d_{\bg}\eta)(g-\bg)\ = \ 0.
\end{equation}
This will yield the proof of Theorem~\ref{theo:nu}. To see that the 
invariant is indeed an invariant of the CR structure, we begin to
check (as the reader can easily convince himself by a straightforward 
computation) that the boundary integral (\ref{eq:intbord}) for $\bg$ 
only depends on the formally determined terms in the asymptotic expansion 
of $\bg$ as the formally undetermined terms are always $O(\exp^{-3r})$ and 
thus do not contribute at infinity. Moreover it is well-known that both
the formally determined part of the Kähler-Einstein metric $\bg$ and
the complex structure in $X× ]R,+\infty)$ are invariants of the CR
structure at infinity only. The desired property is then proved.

\medskip

The computations of the limit in Equation (\ref{eq:limite})
will be broken into two parts. We consider first the case of 
the $\eta$-invariant. 
 
\begin{lem}\label{lem:eta} One has 
$\lim\limits_{r\to\infty}\ (d_{\bg}\eta)(g-\bg)= 0$.
\end{lem}
 
\begin{proof} The first variation of the $\eta$-invariant is the integral
of a local quantity given by the scalar product of the metric variation 
$g-\bg$ against a quadratic form $t^{\bg}$ depending on the third 
derivatives of the 
metric~$\bg$~\cite{AtPatSiIII76}:
\begin{equation}\label{eq:deta}
(d_{\bg}\eta)(g-\bg) = \int_{S_r} \la t^{\bg},g-\bg\ra \,
d\!\vol_{(S_r,\bg)}.
\end{equation}
However, since $\bg$ is asymptotically complex hyperbolic and $g-\bg$
is of order $O(\exp^{-2r})$, the only terms that will contribute at
infinity are $k\exp^{-r}$ ($k$ as in section~\ref{sec:da}; recall this
term is $O(\exp^{-2r})$) and the
highest order terms of $\bg$ and $d\!\vol_{(S_r,\bg)}$. As explained
in Remark~\ref{rem:reexpressed}, these terms have exactly the same
coefficients in a basis adapted to the CR-structure as they would
have in the complex hyperbolic space $\CM{\mathbf H}^2$. Furthermore,
the distance spheres are homogeneous under $U(2)$ in the complex
hyperbolic space, hence is any curvature quantity as 
$t^{\CM{\mathbf H}^2}$. This implies that the restriction of $t^{\bg}$
to the contact distribution $H$ at each point of $S_r$ should be (at 
first order) a multiple of $\gamma$. Its scalar product with $k$ then
vanishes and the lemma is proved. \end{proof}

\medskip

We now manage the local terms in formula (\ref{eq:limite}). 

\begin{lem}\label{lem:B}
One has $\lim\limits_{r\to\infty} \ \int_{S_r} (d_{\bg}B)(g-\bg) \ = \ 0$. 
\end{lem}

\begin{proof} For sake of simplicity, let $\dg = g-\bg$ and $\dR=d_{\bg}R(\dg)$.
The first step is provided by the (obvious) computation:
\begin{equation}\label{eq:variation}\begin{split}
- 12\pi^2\,(d_{\bg}B)(\dg) = & \int_{S_r} 
(d_{\bg}\tv)(\dg)\, (\II^{\bg}\land\II^{\bg}\land\II^{\bg}) 
- 3\,\int_{S_r} (d_{\bg}\tv)(\dg)\, (\II^{\bg}\land R^{\bg})\\
& + \int_{S_r} \tv((d_{\bg}\II)(\dg)\land\II^{\bg}\land\II^{\bg}) 
  + \textrm{ circ. permut. } \\
& - 3\,\int_{S_r} \tv((d_{\bg}\II)(\dg)\land R^{\bg}) 
- 3\,\int_{S_r} \csum ((d_{\bg}\II)(\dg)
(\cdot,R_{.,.}^{\bg}\overline{\mathrm{n}}) \\
& - 3\,\int_{S_r} \tv(\II^{\bg}\land \dR)
- 3\,\int_{S_r} \csum (\II (\cdot,\dR_{.,.}
\overline{\mathrm{n}}) )
\end{split}\end{equation}
(note that the unit normals $\overline{\mathrm{n}}$ to the spheres $S_r$
do not change when passing from $\bg$ to $g$).
As in the proof of Lemma~\ref{lem:eta}, 
decay considerations show that it is enough to treat every 
occurrence of $\dg$ as $k\exp^{-r}$ and to compute every term in $\bg$ 
at highest order. Equivalently, one can consider that $\bg$ is the
standard complex hyperbolic metric and $k$ is a $J_0$-anti-invariant
quadratic form on the contact distribution of the standard
structure. 

Now we perform the calculation using the same notations as in section
\ref{sec:weyl}, using the frame $(e_{i})$ as in \eqref{frame}, the basis 
of 2-forms $(\omega^{i}_{±})$ as in \eqref{omegai}, and the form
\eqref{defk} for $k$.

\smallskip

We shall now show that each of the integrands in formula
(\ref{eq:variation}) 
contributes pointwise as $0$ in the limit. Two of them are handled easily:

{\flushleft (i)} the map 
$\tv$ depends only on the volume form. Since $k$ is tracefree,
the first variation of $\tv$ vanishes.

\smallskip

{\flushleft (ii)} the variation of the second fundamental form 
(up to highest order term) is given by $-\frac{1}{2}\, \exp^{-r} k$. 
An easy computation shows that all the terms involving 
$(d_{\bg}\II)$ contribute as zero (as $\partial_r$ is the normal to 
the geodesic spheres for the model space as well as the modified metric, 
the variation of the second fundamental form only
depends on the behavior of $k$ along the ray $]R,+\infty[×\{ p\}$).

\medskip

It remains to study the terms containing first variations of the
curvature, which deserve slightly more attention.

In the frame $(e_{i})$ given by \eqref{frame}, so that
$\partial_{r}=e_{1}$, one has
$$ \II=e^{2}e_{2}+\frac{1}{2}(e^{3}e_{3}+e^{4}e_{4}).
$$
We may then explicit the terms involving $\dR$ in \eqref{eq:variation}:
noting $v=e^{2}\land e^{3}\land e^{4}$, we obtain
\begin{equation}\label{eqtv}
\begin{split}
\tv(\II \land \dR)
&= v \left\langle e^{2}\land e^{3}\land e^{4},
 e_{2}\land \dR_{e_{3},e_{4}}
 +\frac{1}{2}(e_{3}\land \dR_{e_{4},e_{2}}+e_{4}\dR_{e_{2},e_{3}}) \right\rangle \\
&= v \left(\langle e^{3}\land e^{4},\dR_{e_{3},e_{4}} \rangle 
 +\frac{1}{2}(\langle e^{4}\land e^{2},\dR_{e_{4},e_{2}}\rangle +\langle e^{2}\land
 e^{3},\dR_{e_{2},e_{3}}\rangle )\right),
\end{split}
\end{equation}
where here $\dR$ is seen as a 2-form with values in 2-vectors.
Similarly, one has
\begin{equation}\label{eqcsum}
\begin{split}
\csum (\II (\cdot,\dR_{.,.}e_{1}) )
&= v \left(\II(e_{2},\dR_{e_{3},e_{4}}e_{1})
 +\II(e_{3},\dR_{e_{4},e_{2}}e_{1})
 +\II(e_{4},\dR_{e_{2},e_{3}}e_{1})\right)\\
&= v \left(\langle e_{2},\dR_{e_{3},e_{4}}e_{1} \rangle 
 + \frac{1}{2}(\langle e_{3},\dR_{e_{4},e_{2}}e_{1} \rangle+\langle e_{4},\dR_{e_{2},e_{3}}e_{1} \rangle)\right),
\end{split}
\end{equation}
where now $\dR$ is seen as a 2-form with values in endomorphisms.

As in the previous section, we consider a path of metrics $g_{t}$,
and compute with the help of the symmetric endomorphism 
$u_{t}=(g_{0}^{-1}g_{t})^{\frac{1}{2}}$ which sends the metric
$g_{0}$ to the metric $g_{t}$. Then one has (see for example 
\cite[2.5]{Gau93} or the previous section)
\begin{equation}\label{eqdR}
\dR = \mathcal{R} + \frac{1}{2}\ad_{\dg}R,
\end{equation}
where $\mathcal{R}$ is the derivative of $u_{t}^{-1}R_{t}$.

In our case, we deform the metrics as Einstein metrics, so that 
$\mathcal{R}$
reduces exactly to the variation $w$ of the Weyl tensor. Let us decompose
$w=w^{+}+w^{-}$ with $w^{±}=dW^{±}(\dg)$.
From section \ref{sec:weyl} we know that $w^{+}$ lies in the subbundle 
generated by
$$ (\omega _{+}^{2})^{2}-(\omega _{+}^{3})^{2},\ \omega _{+}^{2}\omega _{+}^{3}+\omega _{+}^{3}\omega _{+}^{2} .
$$
Moreover the variation $w^-$ must be a section of the bundle $\cJ$ 
introduced in the previous section, as the
highest order term of $W^-$ should stay in the kernel bundle of 
$(d^{\nabla})^*d^{\nabla} + d^{\nabla}(d^{\nabla})^*$ when one passes
from $\bg$ to $g$. Hence $w^-$ lives in the bundle generated by
$$
(\omega _{-}^{2})^{2}-(\omega _{-}^{3})^{2},\ \omega _{-}^{2}\omega _{-}^{3}+\omega _{-}^{3}\omega _{-}^{2}.
$$

Now let us understand formula \eqref{eqtv}: in that formula, $\dR$ is
seen as a 2-form with values in 2-vectors, so we have to add to the
variation \eqref{eqdR} the variation of the musical isomorphism
$\Omega^{2}=\Lambda ^{2}TM$.
Actually, as before, it is easy to check that the contribution of the
variation of the musical isomorphism, as well as the contribution of $\ad_{\dg}R$,
vanish. Therefore we are reduced to check the vanishing of
\eqref{eqtv} for the tensor $w$, that is of the quantity
$$
\langle e^{3}\land e^{4},w_{e_{3},e_{4}}\rangle +\frac{1}{2}(\langle e^{4}\land e^{2},w_{e_{4},e_{2}}\rangle +\langle e^{2}\land
 e^{3},w_{e_{2},e_{3}}\rangle ) ;
$$
since $w$ does not involve $\omega^{1}_{±}$, one has
$w_{e_{3},e_{4}}=0$, and we are reduced to study
\begin{align*}
\langle e^{4}\land e^{2},w_{e_{4},e_{2}}\rangle +\langle e^{2}\land e^{3},w_{e_{2},e_{3}}\rangle
&= \frac{1}{4}\left(
\sum _{2}^{3}\langle \omega _{+}^{i},w(\omega _{+}^{i})\rangle +
\sum _{2}^{3}\langle \omega _{-}^{i},w(\omega _{-}^{i})\rangle \right)\\ 
&=0.
\end{align*}

In the same way, we attack formula \eqref{eqcsum}: again, all terms
vanish, except maybe the one induced by $w$, which we are now going to
calculate: using the fact that $w$ is a 2-form with values in
orthogonal endomorphisms, \eqref{eqcsum} becomes
$$
-\langle e_{1},w_{e_{3},e_{4}}e_{2}+\frac{1}{2}
(w_{e_{4},e_{3}}e_{2}+w_{e_{2},e_{3}}e_{4}) \rangle 
$$
which becomes, using the Bianchi identity for $w$,
$$
-\frac{1}{2}\langle e_{1},w_{e_{3},e_{4}}e_{2} \rangle ,
$$
and this vanishes because $w$ does not involve $\omega^{1}_{±}$.
\end{proof}

\bigskip

\section{Relations with the Burns-Epstein invariant}\label{sec:nu-variation}

\bigskip

In the Kähler-Einstein case, on a complex domain, Burns and Epstein 
\cite{BurEps90b} have identified the boundary term of the integral 
\eqref{eq:burns-epstein} as their invariant $\mu$ of the CR boundary,
so that in this case one has
\begin{equation}\label{numu}
 \nu\ =\ 3\mu+2.
\end{equation}

For general CR manifolds with trivial holomorphic bundle (where the 
Burns-Epstein invariant $\mu$ is still defined from \cite{BurEps88})
this relation may not hold, and the difference $\nu-3\mu$
seems difficult to calculate.

We shall give below a first step in this direction.
We will prove that the variation of $\nu$ with respect to any
deformation of the complex structure on the contact
distribution at infinity, is equal to 3 times the variation of the
Burns-Epstein invariant with respect to the same deformation. This
proves that $\nu$ equals $3\mu$ up to a constant whose determination
involves delicate normalization problems, which are out of the scope 
of this paper.

\medskip

\subsection*{Variations of the $\nu$-invariant} We now fix 
a contact structure $H$
on a manifold $X$ of dimension $3$. The set of CR structures compatible
to $H$ is contractible, hence one can always relate any two complex
structures on $H$ by a path. Here we shall study the infinitesimal
variation of $\nu(X)$ when the complex structure varies but the
contact structure remains fixed.

\smallskip

We recall a notation: if $J_0$ is a complex structure on $H$, 
the set of deformations of
$J_0$ is the set of sections of the bundle of anti-$\mathbb{C}$-linear
endomorphisms of $(H,J_0)$, already encountered in section \ref{sec:weyl}.
Any element of this space may be described as a linear map sending the
space of $(0,1)$-vectors (for $J_0$) into the space of $(1,0)$-vectors,
extended on the whole of $TH\otimes\mathbb{C}$ by its
conjugate. As a result, we will sometimes denote any such element $E$ in
$\mathcal{J}$
as $E_{\bar 1}^1$, its expression in any frame $Z_1$ of the $(1,0)$-vectors
in $H$ and the corresponding coframe $\vti$.

\smallskip

\begin{theo}\label{theo:nuderivee}
Let $J_0$ a compatible complex structure on the contact distribution $H$
of a contact $3$-manifold $X$.
Then $\nu(X,H,J)$ is a smooth function in $J$ around $J_0$. Moreover, if
$E$ is a section of $\mathcal{J}$ and $J(t)$ is a curve of complex
structures on $H$ with $J(0)=J_0$, $J'(0)=E$, then
\begin{equation}\label{eq:nuderivee}
\frac{d}{dt}\left( \nu(X,H,J(t)) \right)_{|t=0} = - \frac{3}{8\pi^2}
\, \int_X  \langle Q , E \rangle  
\end{equation}
where $Q$ is the Cartan tensor of the CR structure defined by $(H,J)$
and $\langle\cdot,\cdot\rangle$ the induced Hermitian scalar product.
\end{theo}

\smallskip

\begin{rem} As the Cheng-Lee relative invariant $\mu(J,J')$ has the
same derivative (up to a factor $3$), this result implies 
Theorem \ref{ctheo}.
\end{rem}

\smallskip

\begin{rem} As the Cartan tensor is CR-covariant, the above expression
depends {\it a priori} on the choice of a contact form only in a change
of scale in the choice of the Hermitian metric on $H$, in the volume 
form $\eta\land d\eta$ and in the Cartan tensor $Q$. When put together,
their behaviors with respect to the choice of $\eta$ exactly cancel, 
thus the integral provides a CR-invariant.
\end{rem}

\smallskip

\begin{proof}
Let $M=[R_0,+\infty)× X$ and for each $t\geq 0$, denote by $J(t)$
and $\bg(t)$ the extended complex structure and Kähler-Einstein metric
(formally determined up to order $2$) defined on $]2R,+\infty)× X$
by section \ref{sec:approx}. If $\hat g$ is a fixed metric on 
$[R_0,2R]× X$, then one may find for each $t\geq 0$ a smooth metric
$g(t)$ on $[R_0,+\infty)× X$ such that 
\[ g(t) = \hat g \ \textrm{ on } \ [R_0,R]× X \ \textrm{ and } \ 
g(t) = \bg(t) \ \textrm{ on } \ ]2R,+\infty)× X .\]
Then, for each $t\geq 0$ and $r>2R$,
\begin{equation}\label{eq:nubar}
\nu^{\bg(t)}(r) \ = \ - \left(\chi-3\tau\right)([R_0,r)× X) \ + \ 
\int_{[R_0,r)× X} \beta(g(t)) \ + \ F(\hat g) 
\end{equation}
where $\beta(g(t))$ is the characteristic polynomial in the curvature
of $g(t)$ corresponding to $\chi - 3\tau$ and $F(\hat g)$ is a fixed
boundary term depending only on the choice of $\hat g$ on $\{R_0\}× X$.
Hence, 
\[ \frac{d}{dt}\left( \nu^{\bg(t)}(r) \right)_{|t=0} = 
\int_{[R_0,r)× X} \frac{d}{dt}\beta(g(t)) .\]
From classical Chern-Weil theory, there exists an (explicitely known)
$3$-form $\alpha$ in the curvature of $g(0)$ and infinitesimal variation
of the Levi-Civita connections of $g(t)$ at $t=0$ such that 
$\frac{d}{dt}\beta(g(t)) = d\alpha$. As $g(t)$ is independent of $t$ on 
$[R_0,R]× X$, this yields:
\[ \frac{d}{dt}\left( \nu^{\bg(t)}(r) \right)_{|t=0} = \int_{\{r\}× X}
\alpha \ \ .\]
As the form $\alpha$ is locally computable from $g(0)$ and $g'(0)$ (and
a finite number of 
derivatives thereof), the convergence as $r$ goes to infinity is uniform and
\begin{equation}\label{eq:derivnu}
\frac{d}{dt}\left( \nu(X,H,J(t)) \right)_{|t=0} = 
\lim_{r\to\infty}\int_{\{r\}× X} \alpha \ \ .
\end{equation}
We now determine the precise form of $\alpha$ which will be suitable for 
our needs. As one takes $r>2R$, $\alpha$ can be computed on $S_r$ from the 
characteristic polynomial for $3c_2 - (c_1)^2$ rather than from $\beta$
(this amounts to restrict polynomials on the Lie algebra $\mathfrak{so}(4)$
to the smaller $\mathfrak{u}(2)$). 

Using Newton's sums $s_i(A) = \tr(A^i)$, this polynomial is known to be
equal to $\frac{1}{8\pi^2}(3s_2 - (s_1)^2)$. Let $\Omega(t)_i^j$ be the 
$\mathfrak{u}(2)$-valued curvature $2$-form of $g(t)$ in a local
frame $(S_0,S_1)$, generating the $(1,0)$-vectors for $J(t)$, obtained 
from Gram-Schmidt orthonormalization of $(\partial_r - ie^{-r}R,\, Z=
X - iJ(t)X)$. Then,
\begin{equation} \label{eq:transgression}
\beta(g(t)) = \frac{1}{8\pi^2} \left( 3\,
\Omega(t)_i^j\land\Omega(t)_j^i - 
\Omega(t)_i^i\land\Omega(t)_j^j \right) \end{equation}
(with the usual summation conventions), and, letting $\Omega = \Omega(0)$,
\begin{equation} \label{eq:transgr2} 
\alpha = \frac{1}{4\pi^2} \left( 3\, \phi_i^j\land
\Omega_j^i - \phi_i^i\land\Omega_j^j \right)\end{equation}
where $\phi$ is the matrix-valued first variation at $t=0$ of the 
Levi-Civita connections.

The work done in sections \ref{sec:asymptotics} and \ref{sec:weyl} shows 
that the curvature of $\bg(0)$ may be decomposed into 3 contributions: a 
dominant term given by the coefficients of the Riemann curvature tensor of 
$\mathbb{C}\mathbf{H}^2$ expressed in the orthonormal basis $(S_0,S_1)$,
a first perturbation at order $\exp^{-2r}$ (originating from the term 
$W_2^-$ in $W^-$ studied in section \ref{sec:weyl}), and further terms
of decay $\exp^{-\frac{5}{2}r}$ at least. In the meanwhile, analogous
conclusions for the Levi-Civita connections, hence for $\phi$, may be
dragged from Lemma \ref{lem-nabla} and Corollary \ref{cor-nablaa} in 
section \ref{sec:asymptotics}. 

We now draw several useful conclusions from these remarks: first of all, 
the contribution in the limit (\ref{eq:derivnu}) of 
the fastest-decay terms in the curvature form is zero. Moreover, as our 
metric is 
Kähler-Einstein, the expression (\ref{eq:transgr2}) can be rewritten as
\begin{equation} \label{eq:transgr3} 
\alpha = \frac{1}{4\pi^2} \left( 3\, \phi_i^j\land
\Omega_j^i - \frac{3}{2}\, i\,\phi_i^i\land \omega \right)
\end{equation}
where $\omega$ here stands for the Kähler form of $\bg(0)$. 

We can now 
study the contributions of the two other terms in the curvature $\Omega$.
First of all, it is expected that the contribution of the highest order
term in the curvature is zero, because the standard CR structure is a
critical point of the Burns-Epstein invariant. As we expect our invariant
to be strongly related to their, it would be no surprise that its derivative
at the standard metric of the complex hyperbolic space is zero. This is
indeed easily checked: the action of curvature of the complex hyperbolic 
space on $(1,0)$-vectors is best described as
\[ R^{\mathbb{C}\mathbf{H}^2}_{\xi,\eta} = - M_{\xi,\eta} 
+ \frac{1}{2}\omega\otimes  i\, Id \]
where the endomorphism-valued $2$-form $M$ is
\[ Z \longmapsto M_{\xi,\eta} (Z) = \frac{1}{2}\left( 
g_{\mathbb{C}\mathbf{H}^2}(\eta,Z) \xi^{1,0} - 
g_{\mathbb{C}\mathbf{H}^2}(\xi,Z) \eta^{1,0} \right). \]
Injecting into (\ref{eq:transgr3}), it remains
\begin{equation} \label{eq:transgr4} 
\alpha = \frac{1}{4\pi^2} \left( - 3\, \phi_i^j\land M_j^i + 
3\, \phi_i^j\land \left(\Omega^{(2)}\right)_j^i \right)
+ o(\exp^{-2r}) \end{equation}
where $\Omega^{(2)}$ denotes the difference between the curvature $\Omega$
and the model curvature $2$-form $\Omega^{\mathbb{C}\mathbf{H}^2}$ of the 
complex hyperbolic space, and $M_i^j$ is
the endomorphism-valued $2$-form $M$ seen in matrix form.

We now show that $\phi_i^j\land M_j^i$ yields an exact form, thus
contributing as zero when integrating $\alpha$. Using a local coframe 
$(s^0,s^1)$ of $(1,0)$-forms dual to the orthonormal frame $(S_0,S_1)$, 
one gets
\[ \phi_i^j\land M_j^i\ = \ 
\frac{1}{2}\left( \phi_0^0\land s^0\land s^{\bar 0}
+ \phi_1^1\land s^1\land s^{\bar 1} + \phi_1^0\land s^1\land s^{\bar 0}
+ \phi_0^1\land s^0\land s^{\bar 1} \right) .\]
This may be transformed as follows. For any metric in the
family $\{\bg(t)\}_{t\geq 0}$, one has
\begin{equation}\label{eq:leviciv} \begin{cases}
ds^0 = \Phi^0_0\land s^0 + \Phi^0_1\land s^1 \\
ds^1 = \Phi^1_0\land s^0 + \Phi^1_1\land s^1 
\end{cases}  \end{equation}
where $\Phi$ is the antihermitian matrix-valued Levi-Civita connection 
$1$-form of $\bg(t)$. Notice then that $\phi_i^j = (\Phi_i^j)'(0)$. Hence,
taking the derivative of Equations (\ref{eq:leviciv}) with respect to 
$t$ at $t=0$ yields expressions
for any $\phi_a^b\land s^a$ in terms of $\dot{s}^a$ and 
$\Phi_a^b\land \dot{s}^b$, where $\dot{s}^c$ denotes $(s^c)'(0)$. 
Easy computations using (\ref{eq:leviciv}) 
a second time for barred indices (and the antihermitian character of
$\Phi$) lend eventually to
\begin{equation*} -2\,\phi_i^j\land M_j^i = 
 \dot{s}^0 \land s^{\bar 0}  +  \dot{s}^1 \land s^{\bar 1}
- ds^{\bar 0}\land \dot{s}^0 -  ds^{\bar 1}\land \dot{s}^1 
=  d\left( \dot{s}^0 \land s^{\bar 0}
+ \dot{s}^1 \land s^{\bar 1}\right)
\end{equation*}
which is the expected exact term---the reader may check that this is a
globally defined $3$-form by chasing its variation under a frame change
from $Z$ to $\exp^{iu}Z$. Thus,
\begin{equation} \label{eq:transgr5} 
\alpha = \frac{3}{4\pi^2} \,
 \phi_i^j\land \left(\Omega^{(2)}\right)_j^i \ 
+ \ \textrm{ (exact terms) } \ + \ o(\exp^{-2r}) .\end{equation}
It then remains to study the contribution of the dominant term in 
$\Omega^{(2)}$, which we know to be the highest-order term in $W^-$.
As this term is $O(\exp^{-2r})$ and must be evaluated against $\phi$, only 
the zeroth-order terms in $\phi$ will contribute. Section 
\ref{sec:asymptotics} yields
\begin{eqnarray}
\phi \ = \ \frac{i}{\sqrt{2}}\, \left( \begin{array}{cccc}
0 & - E_{1}^{\bar 1}\theta^1  \\
E_{\bar 1}^{1}\theta^{\bar 1} & q_1^1 
\end{array}\right) \ + \ o(1) 
\end{eqnarray} 
where the first row (column) is the $\{\partial_r,R\}$-complex line and the
second one corresponds the $(1,0)$-part of $H$, generated by an orthonormal 
frame $Z_1$ with associated coframe $\theta^1$, and $q_1^1$ is the first 
variation of the Webster connection induced by the variation of $J$ at 
infinity. From the definition of $\mathcal{J}$, only the non-diagonal
terms have a non-zero contribution against $W_2^-$. Since $W_2^- = 
a\, Q$ (where $Q$ is the Cartan tensor, see section \ref{sec:weyl}), 
this finally proves \begin{equation} \label{eq:transgr6} 
\lim_{r\to\infty} \, \int_{\{r\}× X} \alpha = const. \,
 \int_X \langle Q,E \rangle  .\end{equation}
The constant can be set by comparing with the case of a domain in
$\mathbb{C}^2$. As every computation done above is local, one should
not be able to distinguish this special 
case from the more general (ACHE) one.
The previous work of Burns and Epstein \cite{BurEps90b} implies 
that variations of our invariant equal three times the variations
of the Burns-Epstein invariant for domains, and the variations of the
latter one with respect to the complex structure have been computed by
in Cheng and Lee \cite{CL90}, from which the constant may be borrowed. 
This enables us to conclude the proof of Theorem 
\ref{theo:nuderivee}. \end{proof}

\medskip

\begin{cor}
For each contact structure $H$ on $X$, there is a constant $a(H)$ such that
$\nu(X,H,J) \ = \ 3\,\mu(X,H,J) \ + \ a(H)$.
\end{cor}

\medskip
 
We have proved above that $\nu$ equals $3\mu$ up to an unknown constant
in each component of the set of contact structures on all possible 
$3$-dimensional oriented manifolds. Giving a better result might be 
a difficult task.
To give an idea of the problem, one can describe our approach as a CR
analogue of what the Atiyah-Patodi-Singer index theorem implies for
the determination of the $\eta$-invariant, whereas Burns and Epstein's
is linked to Chern-Simons theory. As is well-known, the
$\eta$-invariant and its Chern-Simons counterpart (lifted from
$\mathbb{R}/\mathbb{Z}$ to $\mathbb{R}$) differ by
normalization constants whose determination is unclear. The constant
appearing in the difference between $\nu$ and $\mu$ may be evidence of 
an analogous phenomenon.
 
\smallskip

\begin{small}{\flushleft\it Acknowledgments}.
The first author thanks M.~Rumin for explaining him the subtleties of
the Burns-Epstein invariant, the second author thanks G.~Lebeau for
his decisive comments, and both are grateful to G.~Carron,
V.~Kharlamov and T.R.~Ramadas for some useful discussions about this work.
\end{small}

\medskip

\end{document}